\documentclass[11pt]{article}
\usepackage{times}
\usepackage{latexsym}
\usepackage{graphicx}
\usepackage{color}
\usepackage{enumitem}
\usepackage{amssymb}

\usepackage[colorlinks=true,linkcolor=magenta,citecolor=blue,pagebackref]{hyperref}

\newcommand{\given}{\mbox{{\bf $\ \mid \ $}}}

\newcommand{\ls}[1]
  {\dimen0=\fontdimen6\the\font \lineskip=#1\dimen0
  \advance\lineskip.5\fontdimen5\the\font \advance\lineskip-\dimen0
  \lineskiplimit=.9\lineskip \baselineskip=\lineskip
  \advance\baselineskip\dimen0 \normallineskip\lineskip
  \normallineskiplimit\lineskiplimit \normalbaselineskip\baselineskip
  \ignorespaces }

\newcommand{\Pf}{\paragraph{{\bf Proof.}}}       
\newcommand{\blot}{\hfill{\vrule height .9ex width .8ex depth -.1ex }}
\newcommand{\EndPf}{\hfill $\blot$ \medskip}     

\newcommand{\argmin}{\mathrm{argmin}}

\iftrue 
\usepackage{amsmath}
\usepackage{amsfonts}
\newcommand{\slim} {\mathop{\rm lim\,sup}}
\newcommand{\ilim} {\mathop{\rm lim\,inf}}
\newcommand{\field}[1]{\mathbb{#1}}
\DeclareMathOperator{\PR}{\field{P}}             
\DeclareMathOperator{\E}{\field{E}}              
\def\N{\field{N}}                                
\def\R{\field{R}}                                
\def\F{\field{F}}                                
\def\A{\field{A}}                                
\def\X{\field{X}}                                
\def\S{{\boldsymbol S}}                                    
\def\Gr{\rm Gr}
\def\Z{\field{Z}}
\def\K{\field{K}}

\def\Sp{\field{S}}
\def\c{\bar c}
\def\n{t}
\def\1-{1-}
\else
\def\PR{\mathop{\rm I\kern -0.20em P}\nolimits}  
\def\E{\mathop{\rm I\kern -0.20em E}\nolimits}   
\def\N{\mathop{\rm I\kern -0.20em N}\nolimits}   
\def\R{\mathop{\rm I\kern -0.20em R}\nolimits}   
\def\F{\mathop{\rm I\kern -0.20em F}\nolimits}   
\def\A{{A}}                                      
\def\X{{X}}                                      
\def\S{{\boldsymbol S}}                                   

\fi

\newtheorem{thm}{Theorem}[section]
\newtheorem{cor}[thm]{Corollary}
\newtheorem{lem}[thm]{Lemma}
\newtheorem{prop}[thm]{Proposition}

\newtheorem{defn}[thm]{Definition}
\newtheorem{rem}[thm]{Remark}
\numberwithin{equation}{section}

\title{On the Convergence of Optimal Actions for Markov Decision Processes and the Optimality of $(s,S)$ Inventory Policies}

\author{Eugene A. Feinberg \\
{\small \it Department of Applied Mathematics and Statistics}\\
{\small \it Stony Brook University}\\ {\small
\it Stony Brook, NY 11794-3600} \\{\small \it
efeinberg@notes.cc.sunysb.edu}
\and
Mark E. Lewis \\{\small \it School of Operations Research and Information Engineering}\\{\small \it Cornell University}\\
{\small \it 221 Rhodes Hall, Ithaca, NY 14853-3801}\\{\small \it mark.lewis@cornell.edu}}

\date{\small \today}

\begin{document}
\maketitle



 \textit{OR/MS subject classifications:} Primary:
 Dynamic Programming/optimal control/Markov/Infinite state; Secondary:
 Inventory/production/Uncertainty/Stochastic

\begin{abstract}
This paper studies  convergence properties of optimal values and actions for
discounted and average-cost Markov Decision Processes (MDPs) with weakly
continuous transition probabilities and applies these properties to the stochastic
periodic-review inventory control problem with backorders, positive setup costs, and
convex holding/backordering costs. The following results are established for  MDPs
with possibly noncompact action sets and unbounded cost functions: (i) convergence
of value iterations to optimal values for discounted problems with possibly non-zero
terminal costs, (ii) convergence of optimal finite-horizon actions to optimal
infinite-horizon actions for total discounted costs, as the time horizon tends to infinity,
and (iii) convergence  of optimal discount-cost actions to optimal average-cost actions
for infinite-horizon problems, as the discount factor tends to 1.

Being applied to the setup-cost inventory control problem, the general results on MDPs  imply the optimality of $(s,S)$ policies and convergence properties of optimal thresholds.  In particular this paper analyzes the setup-cost inventory control problem without two assumptions often used in the literature: (a)  the demand is either discrete or continuous or (b) the backordering cost is higher than the cost of backordered inventory if the amount of backordered inventory is large.  
\end{abstract}
\section{Introduction}\label{intro}
Since Scarf~\cite{scarf} proved the optimality of $(s,S)$ policies for finite-horizon
problems with continuous demand, there have been significant efforts to extend this result to other models. 
Arthur F. Veinott~\cite{vei66a, vei66b} was one of the
pioneers in this exploration, and he combined a deep understanding of Markov decision
processes with a passion for the study of inventory control. It is a great pleasure to
dedicate this paper to him.

This paper  introduces new results on Markov Decision Processes (MDPs) with
infinite state spaces, weakly continuous transition probabilities, one-step costs that
can be unbounded, and possibly noncompact action sets under the discounted and
average-cost criteria. The results on MDPs are applied to the  stochastic
periodic-review setup-cost inventory control problem.  We show that this problem satisfies
general conditions sufficient for the existence of optimal policies, the validity of the
optimality equations, and the convergence of value iterations.  In particular, these
results are used to show the optimality of $(s,S)$ policies for finite-horizon
problems, and for infinite-horizon problems with the discounted and long-term
average-cost criteria.

Since the 1950s,  inventory control has been one of the major motivations for studying
MDPs.  However, until recently there has been a gap between the available results in the MDP theory
and the results needed to analyze inventory control problems. Even now most work on
inventory control assumes that the demand is either discrete or continuous. Moreover,
the proofs are often problem-specific and do not use general results on MDPs, which often provide additional insight.  For example, Theorem~\ref{th:sSdisc} below states convergence properties of optimal thresholds in addition to the existence of optimal $(s,S)$ policies, and the proof of this theorem is based on Theorem~\ref{t:lima} and  Corollary~\ref{l:averopt} established for MDPs.

With such a long history, the inventory control literature is far too expansive to attempt
a complete literature review. The reader is pointed to the books by
Bensoussan~\cite{ben}, Beyer et al.~\cite{BCST}, Heyman and Sobel~\cite{heysol}, Porteus~\cite{Port},
Simchi-Levi et al.~\cite{SCB}, and Zipkin~\cite{Zip}. Applications of MDPs to
inventory control are also discussed in Bertsekas~\cite{bert00}. In the case of
inventory control, under the average cost criterion the optimality of $(s,S)$ policies
was established by Iglehart~\cite{igle63} and Veinott and Wagner~\cite{vei65} in the
continuous and discrete demand cases,
respectively. As explained in Beyer and Sethi~\cite[p. 526]{bey99} in detail, the analysis in Iglehart~\cite{igle63} assumes the existence of a demand density. The proofs for discrete demand distributions were significantly  simplified by Zheng~\cite{zhe91}. 
Zabel~\cite{za} corrected Scarf's~\cite{scarf} results on finite-horizon inventory control.  Beyer and Sethi~\cite{bey99} described and corrected gaps in the proofs in \cite{igle63, vei65} on infinite-horizon inventory control with long-term average costs.  
Almost all studies of infinite-horizon inventory control deal with either discrete or
continuous demand. In some cases, the choice between the use of discrete and continuous distributions depends on a particular application. There is an important practical reason why many studies use discrete
demand. In operations management practice, the overwhelming majority of information systems record
integer quantities of demand and stock level.
 Without assumptions that the demand is discrete or continuous, the optimality of $(s,S)$
policies for average cost inventory control problems follows from Chen and
Simchi-Levi~\cite{che04a}, where under some technical assumptions coordinated
price-inventory control is studied and methods specific to inventory control are used.
Huh et al.~\cite{huh} developed additional problem-specific methods for inventory
control problems with arbitrary distributed demands. Under some additional
assumptions, including the assumption that holding costs are bounded above by
polynomial functions, the optimality of $(s,S)$ policies also takes place when the
demand evolves according to a Markov chain; see Beyer et al.~\cite{BCST} and the
references therein.

Early studies of MDPs dealt with finite-state problems and infinite-state problems
with bounded costs.  The case of average costs per unit time is more difficult than
the case of total discounted costs. Sennott~\cite{sen99} developed the theory for
the average-cost criterion for countable-state problems with unbounded costs.
Sch\"al \cite{sch75, sch93} developed  the theory for uncountable state problems
with  discounted and average-cost criteria when action sets are compact.  In
particular, Sch\"al \cite{sch75, sch93} identified two groups of assumptions on
transition probabilities: weak continuity and setwise continuity. As explained in
Feinberg and Lewis~\cite[Section 4]{FL}, models with weakly continuous transition
probabilities are more natural for inventory control than models with setwise
continuous transition probabilities.
 Hern\'andez-Lerma and
Lasserre~\cite{HL:96} developed the theory for problems with setwise continuous
transition probabilities, unbounded costs, and possibly noncompact action sets.
Luque-Vasques and Hern\'andez-Lerma~\cite{LVHL} identified an additional
technical difficulty in dealing with problems with weakly continuous transition
probabilities even for finite-horizon problems by demonstrating that Berge's
theorem, that ensures semi-continuity of the value function,  does not hold for
problems with noncompact action sets.  Feinberg and Lewis~\cite{FL} investigated
 total discounted costs for inf-compact cost functions and obtained sufficient optimality
conditions for average costs.  Compared to Sch\"al~\cite{sch93} these  results
required an additional local boundedness assumption  that  holds for inventory
control problems, but its verification is not easy. Feinberg et al.~\cite{FKZMDP,
FKZ} introduced a natural class of $\K$-inf-compact cost functions, extended
Berge's theorem to noncompact action sets, and developed the theory of MDPs with
weakly continuous transition probabilities, unbounded costs, and with the criteria of
total discounted costs and long-term average costs. In particular, the results from
\cite{FKZMDP} do not require the validity of the local boundedness assumption.
This simplifies their applications to inventory control problems. Such applications
are considered in Section~\ref{S6} below. The tutorial by Feinberg~\cite{fe16}
describes in detail the applicability of recent results on MDPs to inventory control.

Section \ref{S2} of this paper describes an MDP with an infinite state space, weakly
continuous transition probabilities, possibly unbounded one-step costs, and possibly
noncompact action sets.  Sections~\ref{S3} and \ref{S4} provide the results for
discounted and average cost criteria.  In particular, new results are provided on the
following topics:  (i) convergence of value iterations for discounted problems with
possibly non-zero terminal values (Corollary~\ref{prop:dcoeC}), (ii) convergence of
optimal finite-horizon actions to optimal infinite-horizon actions for total discounted
costs, as the time horizon tends to infinity (Theorem~\ref{t:lima}), and (iii)
convergence  of optimal discount-cost actions to optimal average-cost actions for
infinite-horizon problems, as the discount factor tends to 1
(Theorems~\ref{l:averoptt} and \ref{l:averopttt}).
%
%
%
%
Studying the convergence of value iterations and optimal actions for discounted
costs with non-zero terminal values in this paper is motivated by inventory control.
As was understood by Veinott and Wagner~\cite{vei65}, without additional
assumptions $(s,S)$ policies may not be optimal for problems with discounted
costs, but they are optimal for large values of discount factors.  Even for large
discount factors, $(s,S)$ policies may not be optimal for finite-horizon problems
with discounted cost criteria and zero terminal costs.  However, $(s,S)$ policies are
optimal for such problems with the appropriately chosen nonzero terminal costs, and
this observation is useful for proving the optimality of $(s,S)$ policies for
infinite-horizon problems.

Section~\ref{S5} relates MDPs to problems whose dynamics are defined by
stochastic equations, as this takes place for inventory control.  Section~\ref{S6}
studies the inventory control problem with backorders, setup costs, linear ordering
costs, and convex holding costs and provides two results on the existence of
discounted and average-cost optimal $(s,S)$ policies.  The first result,
Theorem~\ref{th:sSdisc}, states the existence of optimal $(s,S)$ policies for large
discount factors and average costs.  It does not use any additional assumptions,
and the proof is based on adding  terminal costs to finite-horizon problems.  The
second result, Theorem~\ref{th:sSdisc1}, states the existence of optimal $(s,S)$
policies for all discount factors under an additional assumption that it is expensive
to keep  a large backordered  amount of inventory. Such assumptions are often used
in the literature; see  Bertsekas~\cite{bert00}, Beyer et al.~\cite{BCST}, Chen and
Simchi-Levi~\cite{che04b, che04a},   Huh et al.~\cite{huh}, and   Veinott and
Wagner ~\cite{vei65}. Theorems~\ref{th:sSdisc} and \ref{th:sSdisc1} also
describe the convergence properties of optimal thresholds for the following two
cases: (i) the horizon length tends to infinity, and (ii) the discount factor tends to
1.

In the conclusion of the introduction, we would like to mention that the results on MDPs
with weakly continuous transition probabilities, non-compact action sets and
unbounded costs presented in this paper are broadly applicable to a wide range of
engineering and managerial problems. Potential applications include resource
allocation problems, control of workload in queues, and a large variety of inventory
control problems. In particular, the presented results should be applicable to combined pricing-inventory
control and to supply chain management; see \cite{che04b, che04a, SCB}. Moreover,
as mentioned above, the results for MDPs presented below significantly simplify the
analysis of the stochastic cash balance problem investigated in \cite{FL} because the
current results do not require verifying the local boundedness assumption introduced in
\cite{FL}. Instead Theorem~\ref{teor3} below can be employed. The periodic-review
setup-cost inventory control problem was selected as an application in this paper
mainly because it is probably the most highly studied inventory control model. We
provide new results for this classic problem.

\section{Definition of MDPs with Borel State and Action Sets}\label{S2}
Consider a discrete-time Markov decision process with the state
space $\X,$  action space $\A$, one-step costs $c,$ and transition probabilities $q.$ The state space $\X$ and action space $\A$ are both assumed to be
Borel subsets of Polish (complete separable metric) spaces.
 If an action $a\in \A$ is selected at a state $x\in\X,$ then a cost
$c(x,a)$ is incurred, where $c:\X\times\A\to\overline\R=\R\cup\{+\infty\},$ and the
system moves to the next state according to the probability distribution
$q(\cdot|x,a)$ on $\X.$  The function $c$ is assumed to be bounded below and
Borel measurable, and $q$ is a transition probability, that is, $q(B|x,a)$ is a Borel
function on $\X\times\A$ for each Borel subset $B$ of $\X,$ and $q(\cdot|x,a)$ is a
probability measure on the Borel $\sigma$-field of $\X$  or each
$(x,a)\in\X\times\A.$

The decision process proceeds as follows: at time $t=0,1,\ldots$ the current state
of the system, $x_\n$, is observed. A decision-maker decides which action, $a$, to
choose, the cost $c(x,a)$ is accrued, the system moves to the next state  according
to $q(\cdot\given x,a),$ and the process continues. Let
$H_\n=(\X\times\A)^{\n}\times\X$ be the set of histories for $\n=0,1,\ldots\ .$ A
(randomized) decision rule at epoch $t=0,1,\ldots$ is a regular transition probability
$\pi_\n$
from $H_\n$ to $\A.$ 
In other words, (i)
$\pi_\n(\cdot|h_\n)$ is a probability distribution on $\A,$
where $h_\n=(x_0,a_0,x_1,\ldots,a_{\n-1}, x_\n)$
and (ii) for any measurable subset $B \subseteq \A$, the function
$\pi_\n(B|\cdot)$ is measurable on $H_\n.$ A policy $\pi$ is a
sequence $(\pi_0,\pi_1,\ldots)$ of decision rules. Moreover, $\pi$
is called non-randomized if each probability measure
$\pi_\n(\cdot|h_\n)$ is concentrated at one point. A non-randomized
policy is called Markov if all decisions depend only on the current
state and time.  A Markov policy is called stationary if all
decisions depend only on the current state. Thus, a Markov policy
$\phi$ is defined by a sequence $\phi_0,\phi_1,\ldots$ of measurable
mappings $\phi_\n:\X \rightarrow \A.$ 
A stationary policy $\phi$ is defined by a
measurable mapping $\phi:\X \rightarrow \A.$ 

The Ionescu Tulcea theorem (see \cite[p. 140-141]{besh96} or  \cite[p.
178]{HL:96}) implies that an initial state $x$ and a policy
$\pi$ define a unique probability distribution $\PR_x^\pi$ on the
set of all trajectories $H_\infty=(\X\times \A)^\infty$ endowed with
the product $\sigma$-field defined by Borel $\sigma$-fields of $\X$
and $\A.$ Let $\E_x^\pi$ be the expectation with respect to this
distribution. For a finite horizon $N=0,1,\ldots$ and a bounded below measurable function
${\bf F}:\X\to{\overline \R}$ called the terminal value, define the
expected total discounted costs
\begin{align}
    v^{\pi}_{N, {\bf F},\alpha}(x) & :=  \E^\pi_x \left[ \sum_{\n=0}^{N-1}
    \alpha^\n c(x_\n,a_\n)+\alpha^N{\bf F}(x_N)\right], \label{discountpi}
\end{align}
where $\alpha \in [0,1)$, $v_{0,{\bf F},\alpha}^\pi(x)={\bf F}(x),$  $x\in\X.$  When ${\bf F}(x)=0$ for all $x\in\X,$ we shall write $v^\pi_{N,\alpha}(x)$ instead of
$v^\pi_{N,{\bf F},\alpha}(x).$  When
$N=\infty$ and ${\bf F}(x)=0$ for all $x\in\X$, \eqref{discountpi} defines the infinite horizon expected
total discounted cost of $\pi$ denoted by  $v_\alpha^\pi(x)$ instead of $v_{\infty,\alpha}^\pi(x).$ The
average costs per unit time are defined as
\begin{align}
    w^{\pi}(x) & := \limsup_{N \rightarrow \infty} \frac{1}{N} \E^\pi_x \sum_{\n=0}^{N-1}
    c(x_\n,a_\n). \label{averagepi}
\end{align}
For each function $V^\pi(x)=v_{N,{\bf F},\alpha}^\pi(x)$, $v_{N,\alpha}^\pi(x)$,
$v_\alpha^\pi(x)$, or $w^{\pi}(x)$, define the optimal cost

\begin{align}
V(x) & := \inf_{\pi \in \Pi} V^{\pi}(x), \label{optdisc}
\end{align}
where $\Pi$ is the set of all policies. A policy $\pi$ is called
\textit{optimal} for the respective criterion if $V^\pi(x)=V(x)$ for
all $x\in \X$.

We remark that the definition of an MDP usually includes the sets of available
actions $A(x)\subseteq\A,$ $x\in\X.$  We do not do this explicitly because we allow
$c(x,a)$ to be equal to $+\infty.$  In other words, a feasible pair $(x,a)$ is modeled
as a pair with finite costs.  To transform this model to one with feasible
action sets, it is sufficient to consider the sets of available actions $A(x)$ such that
$A(x)\supseteq A_c(x),$ where $A_c(x)=\{a\in\A:c(x,a)<+\infty\},$ $x\in\X.$  In particular, it is possible to set $A(x):=A_c(x),$ $x\in\X.$ In
order to transform an MDP with  action sets $A(x)$ to a MDP with action sets $\A,$
$x\in\X,$ it is sufficient to set $c(x,a)=+\infty$ when $a\in\A\setminus A(x).$ Of
course, certain measurability conditions should hold, but this is not an issue when
the function $c$ is measurable.  We remark that early works on MDPs by
Blackwell~\cite{Black} and Strauch~\cite{Str} considered models with $A(x)=\A$
for all $x\in\X.$ This approach caused some problems with the generality of the
results because the boundedness of the cost function $c$ was assumed and
therefore $c(x,a)\in\R$ for all $(x,a).$  If the cost function is allowed to take
infinitely large values, models with $A(x)=\A$  are as general as models with
$A(x)\subseteq \A,$  $x\in\X.$

\section{Optimality Results for Discounted Cost MDPs with Borel State and Action Sets}\label{S3}

It is well-known (see e.g. \cite[Proposition 8.2]{besh96}) that
$v_{\n,{\bf F}, \alpha}(x)$ satisfies the following \textit{optimality
equation},
\begin{align}\label{fdcoe}
v_{\n+1,{\bf F},\alpha}(x) &= \inf_{a\in \A} \{c(x,a)+\alpha \int
v_{\n,{\bf F},\alpha}(y)q(dy|x,a)\}, &  x \in \X,\ \n=0,1,\ldots\ .
\end{align}
In addition, a Markov policy $\phi$, defined at the first $N+1$ steps by the
mappings $\phi_0,\ldots,\phi_{N}$  satisfying the following equations for all $x\in
\X$ and all $\n=0,\ldots,N,$
\begin{align}\label{fedisopt} v_{\n+1,{\bf F},\alpha}(x) &=
c(x,\phi_{N-\n}(x))+\alpha \int
v_{\n,{\bf F},\alpha}(y)q(dy|x,\phi_{N-\n,\alpha}(x)), & x \in \X,
\end{align}
is optimal for the horizon $N+1;$  see e.g. \cite[Lemma 8.7]{besh96}.

It is also well-known (see e.g. \cite[Proposition 9.8]{besh96}) that
$v_{\alpha}(x)$ satisfies the following \textit{discounted cost
optimality equation},
\begin{align}\label{dcoe}
v_{\alpha}(x) &= \inf_{a\in \A} \{c(x,a)+\alpha \int
v_{\alpha}(y)q(dy|x,a)\}, &  x \in \X.
\end{align}
According to \cite[Proposition 9.12]{besh96}, a stationary policy $\phi^\alpha$ is optimal if and only if
\begin{align}\label{edisopt}
v_{\alpha}(x) &= c(x,\phi^\alpha(x))+\alpha \int
v_{\alpha}(y)q(dy|x,\phi^\alpha(x)), & x \in \X.
\end{align}

However, additional conditions on cost functions and transition probabilities are
needed to ensure the existence of optimal policies. Earlier conditions required
compactness of action sets. They were introduced  by Sch\"al~\cite{sch75} and
consisted of two sets of conditions that required either weak or setwise continuity
assumptions.  For setwise continuous transition probabilities, Hernandez-Lerma and
Lasserre~\cite{HL:96} extended these conditions to MDPs with general action sets
and cost functions $c(x,a)$ that are inf-compact in the action variable $a.$ Feinberg
and Lewis~\cite{FL} obtained results for weakly continuous transition probabilities
and inf-compact cost functions.  Feinberg et al.~\cite{FKZMDP} generalized and
unified the results by Sch\"al~\cite{sch75} and Feinberg and Lewis~\cite{FL} for
weakly continuous transition probabilities to more general cost functions by using
the notion of a $\K$-inf-compact function. $\K$-inf-compact functions were
originally introduced in \cite[Assumption {\bf W*}]{FKZMDP} without using the
term $\K$-inf-compact, and formally introduced and studied in Feinberg et al.
\cite{FKV, FKZ}.  As explained in Feinberg and Lewis~\cite[Section 4]{FL}, weak
continuity holds for periodic review inventory control problems. The setwise
continuity assumption may not hold, but it holds for problems with continuous or
discrete demand distributions.  This paper focuses on the essentially more general
case of weakly continuous transition probabilities.

Let $\mathbb{U}$ be a metric space and $U\subseteq\mathbb{U}.$ Consider a
function $f:U\to\overline{\mathbb{R}}. $ For $V\subseteq U$ define the level sets
\begin{equation}\label{def-D}
\mathcal{D}_f(\lambda;V):=\{y\in V \, : \,  f(y)\le
\lambda\},\qquad \lambda\in\R.
\end{equation}
A function  $f:U\to\overline{\mathbb{R}}$ is called \textit{lower semi-continuous}
at a point $y\in U$ if $f(y)\le\liminf_{n\to\infty}f(y^{(n)})$ for every sequence
$\{y^{(n)}\in U\}_{n=1,2,\ldots}$  converging to $y.$  A function
$f:U\to\overline{\mathbb{R}}$ is called lower semi-continuous if it is lower
semi-continuous at each $y\in U.$  A function $f:U\to\overline{\mathbb{R}}$ is
called \textit{inf-compact}  if all the level sets $\mathcal{D}_f(\lambda;U)$ are
compact. Inf-compact functions are lower semi-continuous.
%
%
%
For three sets $U,$ $V,$ and $W$, where $V\subset U,$ and two functions $g:V\to
W$ and $f:U\to W,$ the function $g$ is called the \textit{restriction} of $f$ to $V$
if $g(x)=f(x)$ when $x\in V.$ 

%

\begin{defn}{\rm (cf. Feinberg et al.~\cite{FKZ,FKV}, Feinberg and Kasyanov~\cite{FK})}
 Let $\Sp^{(i)}$ be metric spaces and $S^{(i)}\subseteq \Sp^{(i)}$ be their nonempty Borel subsets, $i=1,2.$
A function $f: S^{(1)}\times S^{(2)}\to \overline{\mathbb{R}}$ is called
$\K$-inf-compact   if, for any nonempty compact subset
 $K$  of $ S^{(1)},$ the restriction of $f$ to $K\times S^{(2)}$ is inf-compact.
\end{defn}

We are mainly interested in applying this definition to the function $f=c,$ where
$c$ is the one-step cost.  In this case, $\X$ and $\A$ are Borel subsets of the
Polish spaces $\Sp^{(1)}$ and $\Sp^{(2)}$  mentioned in the definition of an MDP.
Inventory control applications often deal with $ \Sp^{(1)} = \X$ and $ \Sp^{(2)}
= \A.$ However, other applications are possible.  For example, assumption~\ref{state:inf2} of
Theorem~\ref{t:cplush}  deals with 
$\Sp^{(1)}=\A$ and $\Sp^{(2)}=\X.$

The next proposition, which follows directly from  Feinberg et al.~\cite[Lemma 2.1]{FKZ}, demonstrates that $\K$-inf-compact cost functions are natural generalizations of inf-compact cost functions considered in Feinberg and Lewis~\cite{FL} and lower semi-continuous cost functions considered in the literature on MDPs with compact action sets, see e.g., Sch\"al~\cite{sch75, sch93}.
\begin{prop}\label{Lemma3.2}
    The following two statements hold:
    \begin{enumerate}[label=(\roman*)]
        \item an inf-compact function   $f:\X\times \A\to \overline{\mathbb{R}}$ is $\K$-inf-compact;
        \item if $A:\X\to 2^\A\setminus\{\emptyset\}$ is a compact-valued upper
            semi-continuous set-valued mapping and $f:\X\times\A\to
            \overline{\mathbb{R}}$ is a lower semi-continuous function such that
            $f(x,a)=+\infty$ for $x\in\X$ and for
            $a\in\A\setminus A(x),$ then the function  $f$ is $\K$-inf-compact, where $2^U$ denotes the set of all subsets of a set $U.$ 
    \end{enumerate}
\end{prop}


\begin{defn}
    The transition probability $q$ is called weakly continuous, if
    \begin{align}
        \label{assm:weak}
            \int_\X f(x)q(dx|x^{(n)},a^{(n)})\to
                \int_\X f(x)q(dx|x^{(0)},a^{(0)}),\qquad {\rm as}\ n\to\infty,
    \end{align}
    for every bounded
    continuous function $f:\X\to\R$ and for each sequence
    $\{(x^{(n)},a^{(n)}), n=1,2,\ldots\}$ on $\X\times\A$ converging to
    $(x^{(0)},a^{(0)})\in\X\times\A.$
\end{defn}

\noindent {\bf Assumption~\bf{W*}.} The following conditions hold:
\begin{enumerate}[label=(\roman*)]
    \item the cost function $c$ is bounded below and $\K$-inf-compact; \label{state:Wstar1} 
    \item if $(x^{(0)},a^{(0)})$ is a limit of a  convergent sequence
        $\{(x^{(n)},a^{(n)}), n=1,2, \ldots\}$ of elements of $\X\times\A$ such
        that $c(x^{(n)},a^{(n)})<+\infty$ for all $n=0,1,2,\ldots,$ then the
        sequence $\{q(\cdot|(x^{(n)},a^{(n)})), n=1,2,\ldots\}$ converges weakly
        to $q(\cdot|(x^{(0)},a^{(0)}));$ that is, \eqref{assm:weak} holds for every
        bounded continuous function $f$ on $\X.$ \label{state:Wstar2}
\end{enumerate}

For example, Assumption~{\bf W*}\ref{state:Wstar2} holds if the transition
probability $q(\cdot|x,a)$ is weakly continuous on $\X\times\A.$  The following
theorem describes the structure of optimal policies, continuity properties of value
functions, and convergence of value iteration.

\begin{thm}\label{prop:dcoe}  {\rm (Feinberg et al.~\cite[Theorem 2]{FKZMDP})}
    Suppose Assumption~{\bf W*} holds. 
    For $\n=0,1,\ldots,$  $N=0,1,\ldots,$ and $\alpha\in [0,1),$ the following statements hold:
    \begin{enumerate}[label = (\roman*)]
      \item the functions $\{v_{\n,\alpha}, t\ge 0\}$ and $v_\alpha$ are
          lower semi-continuous on $\X$, and
          $v_{\n,\alpha}(x)\to v_\alpha (x)$ as $\n \to +\infty$ for each $x\in \X;$
          \label{state:finite1}
      \item the value functions $\{v_{\n,\alpha}, t\ge 0\}$  satisfy the optimality equations
        \begin{align}\label{eq433}
            v_{\n+1,\alpha}(x) & =\min\limits_{a\in \A}\left\{c(x,a)+\alpha
            \int_\X v_{\n,\alpha}(y)q(dy|x,a)\right\},\quad x\in
            \mathbb{X},
        \end{align}
        and the nonempty sets $A_{\n,\alpha}(x):=\{a\in
        \A:\,v_{\n+1,\alpha}(x)=c(x,a)+\alpha \int_\X v_{\n,\alpha}(y)q(dy|x,a)
        \}$, $x\in \X$, satisfy the following properties:
        \begin{enumerate}[label = (\alph*)]
                      \item the graph ${{\Gr}}_\X(A_{\n,\alpha})=\{(x,a):\, x\in\X, a\in
                          A_{\n,\alpha}(x)\}$ is a Borel subset of $\X\times
                          \mathbb{A};$
                      \item the following hold: \begin{enumerate}
                      \item[(b1)] if $v_{\n+1,\alpha}(x)=+\infty$, then $A_{\n,\alpha}(x)=\A;$
                       \item [(b2)]  if $v_{\n+1,\alpha}(x)<+\infty$, then $A_{\n,\alpha}(x)$ is
                          compact;\end{enumerate}
        \end{enumerate} \label{state:markov1}
      \item for each horizon $(N+1),$ there exists a Markov optimal 
          policy $(\phi_0,\ldots,\phi_{N});$ \label{state:markov2}
      \item if for an $(N+1)$-horizon Markov policy $(\phi_0,\ldots,\phi_{N})$ the
          inclusions $\phi_{N-\n}(x)\in A_{\n,\alpha}(x)$, $x\in\X,$ $\n=0,\ldots,N$
          hold, then this policy is $(N+1)$-horizon optimal;
      \item  the value function $v_\alpha$  satisfies the optimality equation
        \begin{align}\label{eq5a}
            v_{\alpha}(x) & =\min\limits_{a\in \A}\left\{c(x,a)+\alpha\int_\X
                v_{\alpha}(y)q(dy|x,a)\right\},\qquad x\in \X;
        \end{align}
      \item the nonempty sets $A_{\alpha}(x):=\{a\in
        \A:\,v_{\alpha}(x)=c(x,a)+\alpha\int_\X
                v_{\alpha}(y)q(dy|x,a) \}$, $x\in \X$, satisfy
        the following properties: 
        \begin{enumerate}
            \item the graph ${\Gr}_\X(A_{\alpha})=\{(x,a):\, x\in\X, a\in
                A_\alpha(x)\}$ is a Borel subset of $\X\times \mathbb{A};$
            \item if $v_{\alpha}(x)=+\infty$, then $A_{\alpha}(x)=\A$ and, if
                $v_{\alpha}(x)<+\infty$, then $A_{\alpha}(x)$ is compact;
        \end{enumerate}
      \item for the infinite horizon there exists a stationary discount-optimal policy
          $\phi_\alpha$, and a stationary policy is optimal if and only if
          $\phi_\alpha(x)\in A_\alpha(x)$ for all $x\in \X;$ \label{state:infinite}
      \item {\rm (Feinberg and Lewis~\cite[Proposition 3.1(iv)]{FL})} if the cost
          function $c$ is inf-compact, then  the functions $v_{\n,\alpha}$,
          $\n=1,2,\ldots$, and $v_\alpha$ are inf-compact on $\X$. \label{state:fl1}
    \end{enumerate}
\end{thm}

\noindent The following corollary extends the previous theorem to nonzero terminal values $\bf F.$ This extension is useful for the analysis of inventory control
problems.

\begin{cor}\label{prop:dcoeC}
    Let Assumption~{\bf W*} hold. 
    Consider a bounded below, lower
    semi-continuous function ${\bf F}:\X\to\overline\R$. The following statements hold for $\n=0,1,2,\ldots,$ $N=0,1,2,\ldots,$ and $\alpha\in [0,1):$
    \begin{enumerate}[label = (\roman*)]
      \item  the  functions $v_{\n,{\bf F},\alpha}$ are bounded below and
          lower semi-continuous; \label{cor:finite1}

      \item the value functions $v_{\n+1,{\bf F},\alpha}$ satisfy the optimality equations
        \begin{align}\label{eq43333}
            v_{\n+1,{\bf F},\alpha}(x) & =\min\limits_{a\in \A}\left\{c(x,a)+\alpha
                \int_\X v_{\n,{\bf F},\alpha}(y)q(dy|x,a)\right\},\quad x\in
                \mathbb{X},
        \end{align}
        where $v_{0,{\bf F},\alpha}(x)={\bf F}(x)$ for all $x\in \X;$
      \item the nonempty \label{eqsetaf}
        sets
        \begin{align*}
            A_{\n,{\bf F},\alpha}(x) & :=\{a\in \A:\,v_{\n+1,{\bf F},\alpha}(x)=c(x,a)+\alpha \int_\X
                v_{\n,{\bf F},\alpha}(y)q(dy|x,a) \},\quad  x\in \X, 
        \end{align*}
        satisfy the following properties:
        \begin{enumerate}[label = (\alph*)]
          \item the graph ${{\Gr}}_\X(A_{\n,{\bf F},\alpha})=\{(x,a):\, x\in\X, a\in
              A_{\n,{\bf F},\alpha}(x)\}$ is a Borel subset of $\X\times
              \mathbb{A};$
          \item the following hold: \begin{enumerate}

          \item[(b1)] if $v_{\n+1,{\bf F},\alpha}(x)=+\infty$, then $A_{\n,{\bf
              F},\alpha}(x)=\A$;
              \item[(b2)] if $v_{\n+1,{\bf F},\alpha}(x)<+\infty$, then
              $A_{\n,{\bf F},\alpha}(x)$ is compact; \end{enumerate}
        \end{enumerate} \label{cor:markov1}
      \item for an $(N+1)$-horizon  problem with the terminal value function ${\bf F},$  there exists a Markov optimal optimal policy
          $(\phi_0,\ldots,\phi_{N})$ and if, for an $(N+1)$-horizon Markov policy
          $(\phi_0,\ldots,\phi_{N})$ the inclusions $\phi_{N-\n}(x)\in A_{\n,{\bf
          F},\alpha}(x)$, $x\in\X,$ $\n=0,\ldots,N,$ hold then this policy is
          $(N+1)$-horizon optimal; \label{cor:markov2}
      \item if ${\bf F}(x)\le v_\alpha(x)$ for all $x\in\X,$ then $v_{\n,{\bf F},\alpha}(x)\to
          v_\alpha (x)$ as $n \to +\infty$ for all $x\in \X;$ \label{cor:terminal}
      \item if the cost function $c$ is inf-compact, then each of the functions $v_{\n,{\bf
          F},\alpha},$ $\n=1,2,\ldots,$ is inf-compact. \label{cor:inf-compact}
    \end{enumerate}
\end{cor}
\Pf  Statements \ref{cor:finite1}-\ref{cor:markov2} are corollaries from statements
\ref{state:finite1}-\ref{state:markov2} of Theorem~\ref{prop:dcoe}. Indeed, the
statements of Theorem~\ref{prop:dcoe}, that deal with the finite horizon $N,$ hold
when one-step costs at different time epochs vary. In particular, if the one-step cost at
epoch $\n=0,1,\ldots,N$ is defined by a bounded below, measurable cost function
$c_\n$ rather than by the function  $c.$ This case can be reduced to the single
function $c$ by replacing the state space $\X$ with the state space
$\X\times\{0,1,\ldots,N\},$  setting $c((x,\n),a)=c_\n(x,a),$ and applying the
corresponding statements of Theorem~\ref{prop:dcoe}.  In our case,
$c_\n(x,a)=c(x,a)$ for $\n=0,1,\ldots,N,$ and $c_{N}(x,a)=c(x,a)+\int_\X {\bf
F}(y)q(dy|x,a).$ The function $c_{N}$ is bounded below and lower semi-continuous.

To prove \ref{cor:terminal} and \ref{cor:inf-compact}, consider first the case when
the functions $c$ and ${\bf F}$ are nonnegative.  In this case,
\begin{align}\label{e:estwT}
    v_{\n,\alpha}(x)\le v_{\n,{\bf F},\alpha}(x)\le v_{\n,v_\alpha,\alpha}(x)
        & =v_\alpha(x), \qquad x\in\X,\ \n=0,1,\ldots\ .
\end{align}
Therefore, for nonnegative cost functions, Statement \ref{cor:terminal} follows from
Theorem~\ref{prop:dcoe}\ref{state:finite1}. Statement \ref{cor:inf-compact}
follows from \ref{cor:terminal}, Theorem~\ref{prop:dcoe}\ref{state:fl1}, and the
fact that $v_{\n,{\bf F},\alpha}\ge v_{\n,\alpha}$ since ${\bf F}$ is nonnegative.  In a general
case, consider a finite positive constant $K$ such that the functions $c$ and ${\bf F}$
are bounded below by $(-K).$ If the cost functions $c$ and ${\bf F}$ are increased by
$K,$ then the new cost functions are nonnegative, each finite-horizon value function
$v_{\n,{\bf F},\alpha}$ is increased by the constant
$d_\n=K(1-\alpha^{\n+1})/(1-\alpha),$ and the infinite-horizon value function
$v_{\alpha}$ is increased by the constant $d=K/(1-\alpha).$ Since $d_\n\le d$ and
$d_\n\to d$ as $\n\to\infty,$ the general case follows from the case of non-negative
cost functions. \EndPf

While Theorem~\ref{prop:dcoe} and Corollary~\ref{prop:dcoeC} state the
convergence of value functions and describe the structure of optimal sets of actions,
the following theorem describes convergence properties of optimal actions. For $x\in\X$ and $\alpha\in [0,1),$
 define the sets $D^*_\alpha(x):=\{a\in \A: c(x,a)\le v_{\alpha}(x)\}.$

\begin{thm}\label{t:lima} Let Assumption~{\bf W*} hold and $\alpha\in [0,1).$ 
Suppose ${\bf F}:\X\to\overline \R$ is bounded below, lower semi-continuous, and
such that for all $x\in\X$
\begin{equation}\label{eq:spdouco}
 {\bf F}(x)\le v_\alpha(x)\qquad{\rm and} \qquad v_{1,{\bf F},\alpha}(x)\ge {\bf F}(x).
\end{equation}
For $x\in\X,$ such that $v_\alpha(x)<\infty,$ the following two statements hold:
\begin{enumerate}[label = (\roman*)]
    \item  the set $D^*_\alpha(x)$ is compact, and 
    $A_{\n,{\bf
        F},\alpha}(x)\subseteq D^*_\alpha(x)$ for all $t=1,2,\ldots,$  where the sets
        $A_{\n,{\bf F},\alpha}(x)$ are defined in Corollary~\ref{prop:dcoeC}\ref{eqsetaf}; \label{state:d-alpha}
    \item each sequence $\{a^{(\n)}\in A_{\n,{\bf F},\alpha}(x), \n=1,2,\ldots\}$ is
        bounded, and all its limit points belong to $A_\alpha(x).$ \label{state:limit-points}
\end{enumerate}
\end{thm}

In particular, if $c(x,a)\ge 0$ for all $x\in\X$, $a\in\A,$ then the function ${\bf F}(x)\equiv 0$ satisfies conditions~\eqref{eq:spdouco}.
In order to prove Theorem~\ref{t:lima}, we need the following lemma, which is a simplified version of \cite[Lemma 4.6.6]{HL:96}.
\begin{lem}\label{l:useful}
Let $A$ be a compact subset of $\A$ and $f,$ $f_n:A\to {\overline\R},$
$n=1,2,\ldots,$ be nonnegative, lower semi-continuous, real-valued functions  such
that $f_n(a)\uparrow f(a)$ as $n\to \infty$ for all $a\in A.$ Let
$a^{(n)}\in\argmin_{a\in A} f_n(a),$  $n=1,2,\ldots,$ and $a^*$
be a limit point of the sequence $\{a^{(n)}, n=1,2,\ldots\} .$  Then
$a^*\in\argmin_{a\in A} f(a).$
\end{lem}
\Pf Let $a^\prime\in\argmin_{a\in A} f(a).$ Then
$f(a^\prime)\ge f_n(a^{(n)}) \ge f_k(a^{(n)})$ for all $n\ge k.$   Since $A$ is compact, then $a^*\in\A.$ Lower
semi-continuity of $f$ and the previous inequalities imply $f_k(a^*)\le
\liminf_{n\to\infty}f_n(a^{(n)})\\ \le f(a^\prime).  $  Thus $f(a^\prime)\ge
f_k(a^*)\uparrow f(a^*).$ Since $f(a^*)\le f(a^\prime),$ then $a^*\in\argmin_{a\in A} f(a).$ \EndPf \\

\noindent {\bf Proof of Theorem~\ref{t:lima}.} We assume without loss of
generality that the bounded below functions $c$ and ${\bf F}$ are nonnegative.  We can
do this because of the arguments provided at the end of the proof of
Corollary~\ref{prop:dcoeC} and the additional argument that, if the one-step cost
functions $c$ and terminal cost functions are shifted by constants, then the set of
optimal finite-horizon action $A_{\n,{\bf F},\alpha}(\cdot)$ and infinite-horizon actions $A_\alpha(\cdot)$ remain unchanged. 

Fix $x\in\X.$ 
Since the function $v_{\n,{\bf F},\alpha}$ is nonnegative  and,
in view of \eqref{e:estwT},  $v_{\n+1,{\bf F},\alpha}(x)\le v_{\alpha}(x),$ 


\begin{align*}
    A_{\n,{\bf F},\alpha}(x) & =\{a\in \A: c(x,a)+\alpha\int_\X v_{\n,{\bf F},\alpha}(y)q(dy|x,a)
        = v_{\n+1,{\bf F},\alpha}(x) \} \subseteq D^*_\alpha(x),\qquad   \n=1,2,\ldots\ .
\end{align*}
Statement \ref{state:d-alpha} is proved. Since $D^*_\alpha(x)$ is compact, every sequence  $\{a^{(\n)}\in
A_{\n,{\bf F},\alpha}(x)\}_{\n=1,2,\ldots}$ is bounded and has a limit point. The theorem follows from
Lemma~\ref{l:useful} applied to the set $A:=D^*_\alpha(x)$ and functions
\begin{align*}
f(a)&=c(x,a)+\alpha \int_\X v_\alpha(y)q(dy|x,a),\qquad\quad\quad a\in A,\\
f_\n(a)&=c(x,a)+\alpha\int_\X v_{\n,{\bf F},\alpha}(y)q(dy|x,a),\qquad \ a\in A,\  \n=0,1,\ldots\ .
\end{align*}%

To verify the conditions of Lemma~\ref{l:useful}, observe that for all $z\in\X$
\[v_\alpha(z)=v_{\n,v_\alpha,\alpha}(z)\ge v_{\n,{\bf F},\alpha}(z)\ge v_{\n,\alpha}(z)\uparrow v_\alpha(z),
\]
where the first equality follows from the optimality equation, the first and the second inequalities follow from $v_\alpha(\cdot)\ge {\bf F}(\cdot)\ge 0,$ and the convergence is stated in Theorem~\ref{prop:dcoe}\ref{state:finite1}; this convergence is monotone because $c$ and ${\bf F}$ are nonnegative functions.  The inequality $v_{1,{\bf F},\alpha}(\cdot)\ge {\bf F}(\cdot)$ in \eqref{eq:spdouco}, equality \eqref{eq43333}, and standard induction arguments imply $v_{\n+1,{\bf F},\alpha}(\cdot)\ge v_{\n,{\bf F},\alpha}(\cdot),$ $\n=0,1,\ldots\ .$ 
Thus Assumption \eqref{eq:spdouco} 
implies that $v_{\n,{\bf F},\alpha}\uparrow v_\alpha,$ and the monotone
convergence theorem implies $f_\n\uparrow f$ as $\n\to\infty.$  \EndPf

\section{Average-Cost MDPs with Borel State  and Action Sets}\label{S4}

The average cost case is more subtle than the case of expected total discounted
costs. The following assumption was introduced by Sch\"al~\cite{sch93}.  Without
this assumption the problem is trivial because $w(x)=+\infty$ for all $x\in\X,$ and
therefore every policy is optimal.\\

\noindent \textbf{Assumption {\bf G}}. $ w^*:=\inf\limits_{x\in
\X}w(x)<+\infty$.\\

Assumption {\bf G} is equivalent to the existence of $x\in\X$ and $\pi\in\Pi$ with
$w^\pi(x)<\infty.$   Define the following quantities for $\alpha\in [0,1)$:
\begin{align*}
    m_{\alpha} & =\inf\limits_{x\in \X}v_{\alpha}(x),\quad
        u_{\alpha}(x)=v_{\alpha}(x)-m_{\alpha},
\end{align*}
\begin{align*}
    \underline{w} & =\ilim\limits_{\alpha\uparrow
        1}(1-\alpha)m_{\alpha},\quad\overline{w}=\slim\limits_{\alpha\uparrow
        1}(1-\alpha)m_{\alpha}.
\end{align*}
Observe that $u_\alpha(x)\ge 0$ for all $x\in\X.$ According to
Sch\"al \cite[Lemma 1.2]{sch93}, Assumption~{\bf G} implies
\begin{align}\label{eq:sch93}
    0\le \underline{w}\le \overline{w}\le w^*< +\infty.
\end{align}
Moreover, Sch\"al \cite[Proposition 1.3]{sch93}, states that, if there exist a
measurable function $u:\ \X\to \R^+,$ where $\R^+:=[0,+\infty),$ and a stationary
policy $\phi$ such that
\begin{align}\label{acoi}
    \underline{w}+u(x) & \geq c(x,\phi(x))+\int u(y)q(dy|x,\phi(x)), & & x\in
        \X,
\end{align}
then $\phi$ is average cost optimal and $w(x)=w^*$ for all $x\in\X$. The following
condition plays an important role for the validity of \eqref{acoi}.\\

\noindent {\bf Assumption {\bf B}.}\ \  Assumption~\textbf{G} holds and
$\sup_{\alpha \in [0,1)} u_\alpha (x) <\infty$ for all $x\in\X.$\\

 We note that the second part of Assumption \textbf{B} is Condition
\textbf{B} in Sch\"al~\cite{sch93}. Thus, under Assumption~\textbf{G}, which is
assumed throughout \cite{sch93}, Assumption \textbf{B} is equivalent to
Condition~\textbf{B} in \cite{sch93}.

For $x\in\X$ and for a nonnegative lower semi-continuous function $u:\X\to \R^+,$ define the set
 \begin{equation}\label{defsetA*}
 A_u^*(x):=\left\{a\in \A\, : \,
\underline{w}+u(x)\ge c(x,a)+\int_\X  u(y)q(dy|x,a) \right\}, \
x\in\X.
\end{equation}
 A stationary policy $\phi$ satisfies (\ref{acoi}) if and only if   $A^*_u(x)\ne\emptyset$  and      $\phi(x)\in  A^*_u(x)$
 for all  $x\in\X.$

Following Feinberg et al.~\cite[Formula (21)]{FKZMDP}, define
\begin{align}\label{eq131}
    u(x) & :=\ilim\limits_{(y,\alpha)\to (x,1-)}u_{\alpha}(y),\quad x\in\X.
\end{align}
In words, ${u}(x)$ is the largest number such that ${u}(x)\le
\liminf_{n\to\infty}u_{\alpha_n}(y_n)$ for all sequences $\{y_n, n\ge 1\}$ and
$\{\alpha_n, n\ge 1\}$ such that  $y_n\to x$ and $\alpha_n\to 1.$

Following Sch\"al~\cite[Page 166]{sch93}, where the notation $\underline w$ is used instead of $\tilde u,$
and  Feinberg et al.~\cite[Formula (38)]{FKZMDP},  for a particular sequence
$\alpha_n\to 1-,$ define
\begin{align}\label{eq1311}
    {\tilde u}(x) & :=\ilim\limits_{(y,n)\to (x,\infty)}u_{\alpha_n}(y),\quad x\in\X.
\end{align}
In words, ${\tilde u}(x)$ is the largest number such that ${\tilde u}(x)\le \liminf_{n\to\infty}u_{\alpha_n}(y_n)$ for all sequences $\{y_n,n\ge 1\}$ converging to $x.$

It follows from these definitions that $u(x)\le{\tilde u}(x),$ $x\in\X.$ However, the
questions, whether $u={\tilde u}$ and whether the values of $\tilde u$ depend on a
particular choice of the sequence $\alpha_n,  $ have not been investigated. If
Assumption~{\bf B} holds, then ${\tilde u}(x)<+\infty, $ $x\in\X.$ If
Assumption~{\bf B} holds and
 the cost function $c$
is inf-compact, then the functions $v_\alpha,$  $u,$ and $\tilde u$ are inf-compact
as well; see Theorem~\ref{prop:dcoe}\ref{state:finite1} for  this fact for $v_\alpha$ and
Feinberg et al.~\cite[Theorem 4(e) and Corollary 2]{FKZMDP} for $u$ and $\tilde
u$.

\begin{thm}\label{teor3} {\rm (Feinberg et al.~\cite[Theorem 4 and Corollary 2]{FKZMDP})}.
    Suppose Assumptions~{\bf W*} and {\bf B} hold. The following two properties hold
    for the function $u$ defined in \eqref{eq131} and for $u=\tilde u,$ where $\tilde u$
    is defined in \eqref{eq1311} for a sequence $\{\alpha_n, n\ge 1\}$ such that $\alpha_n\uparrow 1:$
     \begin{enumerate}[label=(\alph*)]
     \item for each $x\in\X$ the set $A_u^*(x)$ is nonempty and compact;
     \item the graph ${{\Gr}_\X}(A_u^*)=\{(x,a):\, x\in\X, a\in A_u^*(x)\}$ is a
              Borel subset of $\X\times \mathbb{A}.$
       \end{enumerate}
    Furthermore, the following statements hold:
    \begin{enumerate}[label = (\roman*)]
      \item
    there exists a stationary policy $\phi$ satisfying \eqref{acoi};
    \item every  policy $\phi$ satisfying \eqref{acoi} is optimal for the average cost
        per unit time criterion, and

    %
    \begin{align} \label{eq5.16}
        w^\phi(x) & =w(x)=w^*=\underline{w}=\overline{w}=\lim\limits_{\alpha\uparrow
            1}(1-\alpha)v_\alpha(x)=\lim\limits_{N\to\infty}\frac{1}{N} \E^\pi_x \sum_{\n=0}^{N-1}
    c(x_\n,a_\n),\qquad
            x\in\X.
    \end{align}
      \end{enumerate}
\end{thm}

If the one-step cost function $c$ is inf-compact,  the minima of
functions $v_\alpha$ possess additional properties.
Set
\begin{align}\label{e:defxa}
    X_\alpha & :=\{x\in\X\,: v_\alpha(x)=m_\alpha\},\qquad \alpha\in [0,1).
\end{align}
Since $X_\alpha=\{x\in\X\,: v_\alpha(x)\le m_\alpha\},$ this set is closed if Assumptions~{\bf G} and {\bf W*} hold.  If the
function $c$ is inf-compact then inf-compactness of $v_\alpha$ implies that the
sets $X_\alpha$ are nonempty and compact.  The following fact is useful for
verifying the validity of Assumption~{\bf B}; see Feinberg and Lewis~\cite[Lemma
5.1]{FL} and the references therein.

\begin{thm}\label{lemC} (Feinberg et al.~\cite[Theorem 6]{FKZMDP}).
Let Assumptions~{\bf G} and {\bf W*} hold. If the function $c$ is inf-compact, then
there exists a compact set $\mathcal{K}\subseteq\X$ such that
$X_\alpha\subseteq \mathcal{K}$ for all $\alpha\in [0,1).$
\end{thm}

According to Feinberg et al.~\cite[Theorem 5 and Corollary 3]{FKZMDP}, certain
average cost optimal policies can be approximated by discount optimal policies with a
vanishing discount factor.  The following theorem describes particular constructions of
such approximations. Recall that, for the function $u(x)$ defined in \eqref{eq131}, for
each $x\in\X$ there exist sequences $\{\alpha_n ,n \geq 1\}$ and $\{x^{(n)}, n \geq
1\}$ such that $\alpha_n\uparrow 1$ and $x^{(n)}\to x,$ where $x^{(n)}\in\X,$
such that $u(x)=\lim_{n\to\infty} u_{\alpha_n}(x^{(n)} ).$ Similarly, for a sequence
$\{\alpha_n , n\geq 1\}$ such that $\alpha_n\uparrow 1$ consider the function $\tilde
u$ defined in \eqref{eq1311}. For each $x\in\X$ there exist a sequence $\{x^{(n)}, n
\geq 1\}$ of points in $\X$ converging to $x$ and a subsequence $\{\alpha^*_n, n
\geq 1\}$ of the sequence $\{\alpha_n, n \geq 1\}$ such that
$\tilde{u}(x)=\lim_{n\to\infty} u_{\alpha^*_n}(x^{(n)} ).$
\begin{thm}\label{l:averoptt}
    Let Assumptions {\bf W*} and {\bf B} hold. For $x\in\X$ and $a^*\in\A,$ the
    following two statements hold:
    \begin{enumerate}[label = (\roman*)]
      \item Consider a sequence $\{( x^{(n)},\alpha_n), n \geq 1\}$ with $0\le
          \alpha_n\uparrow 1,$ $ x^{(n)}\in\X,$   $ x^{(n)}\to x,$   and
          $u_{\alpha_n}(x^{(n)})\to u(x)$ as $n\to\infty$. If there are  sequences of
          natural numbers $\{n_k ,k \geq 1\}$ and actions
          $\{a^{(n_k)}\in A_{\alpha_{n_k}}(x^{(n_k)}), k \geq 1\},$ such that
          $n_k\to\infty$ and $ a^{(n_k)}\to a^*$ as $k\to\infty,$ then $a^* \in
          A^*_u(x),$ where the function $u$ is defined in \eqref{eq131};
          \label{state:converge1}
      \item Suppose $\{\alpha_n, n \geq 1\}$ is a sequence of discount factors such that
          $\alpha_n\uparrow 1.$ Let $\{\alpha^*_n, n \geq 1\}$ be its subsequence and
          $\{x^{(n)}, n \geq 1\}$ be a sequence of states such that $x^{(n)}\to x$
          and  $u_{\alpha^*_{n}}(x^{(n)})\to \tilde{u}(x)$ as
          $n\to\infty,$ where the function $\tilde u$ is defined in \eqref{eq1311} for
          the sequence $\{\alpha_{n},n \geq 1\}$. If there are actions $a^{(n)}\in
          A_{\alpha^*_{n}}(x^{(n)})$ such that $ a^{(n)}\to a^*$ as $n\to\infty,$
          then $a^*\in A^*_{\tilde u}(x).$ \label{state:converge2}
    \end{enumerate}
\end{thm}
\Pf To show \ref{state:converge1}, consider sequences whose existence is assumed in
the theorem. We have
\begin{align*}
    v_{\alpha_{n_k}}(x^{(n_k)}) & =c(x^{(n_k)},a^{(n_k)})
        +\alpha\int_\X v_{\alpha_{n_k}}(y)q(dy|x^{(n_k)},a^{(n_k)}).
\end{align*}
This implies (with a little algebra)
\begin{align*}
    u_{\alpha_{n_k}}(x^{(n_k)})+(1-\alpha_{n_k})m_{\alpha_{n_k}} & =
        c(x^{(n_k)},a^{(n_k)}) +\alpha\int_\X u_{\alpha_{n_k}}(y)q(dy|x^{(n_k)},a^{(n_k)}).
\end{align*}
Fatou's lemma for weakly converging measures (see e.g., Feinberg et
al.~\cite[Theorem 1.1]{FKZTV}), the choice of the sequence $x^{(n_k)},$ and
Theorem~\ref{teor3} yield
\begin{align*}
    \underline {w}+u(x)\ge c(x,a^*)+ \int_\X
        u(y)q(dy|x,a^*).
\end{align*}
Thus $a^*\in A^*_u(x).$  The proof of Statement \ref{state:converge2} is similar.
\EndPf

\begin{cor}\label{l:averopt}
    Let Assumptions {\bf W*} and {\bf B} hold.  For $x\in\X$ and $a^*\in \A,$ the
    following hold:
    \begin{enumerate}[label = (\roman*)]
      \item if each sequence $\{(\alpha^*_n, x^{(n)}), n \geq 1\}$ with $0\le
          \alpha^*_n\uparrow 1,$ $x^{(n)}\in\X,$ and $ x^{(n)}\to x,$ contains a
          subsequence $(\alpha_{n_k}, x^{(n_k)}),$ such that there exist actions
          $a^{(n_k)}\in A_{\alpha_{n_k}}(x^{(n_k)})$ satisfying $ a^{(n_k)}\to
          a^*$ as $k\to\infty,$ then $a^*\in A^*_u(x)$ with the function $u$ defined
          in \eqref{eq131}; \label{state:averopt1}
      \item if there exists a sequence $\{\alpha_n, n \geq 1\}$ such that
          $\alpha_n\uparrow 1$ and for every sequence of states $\{x_n\to x\}$ from
          $\X$ there are actions $a^{n}\in A_{\alpha_n}(x^{(n)}),$ $n=1,2,\ldots,$
          satisfying $a_n\to a^*$ as $n\to\infty,$ then $a^*\in A^*_{\tilde u}(x),$
          where the function $\tilde u$ is defined in \eqref{eq1311} for the sequence
          $\{\alpha_{n}, n \geq 1\}.$ \label{state:averopt2}
    \end{enumerate}
\end{cor}
\Pf Statement \ref{state:averopt1} follows from
Theorem~\ref{l:averoptt}\ref{state:converge1} applied to a sequence
$\{(\alpha^*_n, x^{(n)}), n \geq 1\}$ with the property $u(x)=\lim_{n\to\infty}
u_{\alpha^*_n}(x^{(n)}).$ Statement \ref{state:averopt2} follows from
Theorem~\ref{l:averoptt}\ref{state:converge2} applied to a sequence $\{ x^{(n)},
n \geq 1\}$  and a subsequence $\{\alpha^*_n, n \geq 1\}$ of $\{\alpha_n, n \geq
1\}$ such that ${\tilde u}(x)=\lim_{n\to\infty} u_{\alpha^*_n}(x^{(n)}).$  \EndPf

The following theorem is useful for proving asymptotic properties of optimal actions for
discounted problems when  the discount factor tends to 1.
\begin{thm}\label{l:averopttt}
    Let Assumptions {\bf W*} and {\bf B} hold. For $x\in\X$  the
    following hold:
    \begin{enumerate}[label = (\roman*)]
    \item There exists a compact set $D^*(x)\subseteq \A$ such that
        $A_\alpha(x)\subseteq D^*(x)$ for all $\alpha\in[0,1);$\label{s:compact}
    \item If $\{\alpha_n, n \geq 1\}$ is a sequence of discount factors $\alpha_n\in
        [0,1),$ then every sequence of infinite-horizon $\alpha_n$-discount cost
        optimal actions $\{a^{(n)}, n \geq 1\},$ where $a^{(n)}\in
        A_{\alpha_n}(x),$ is bounded and therefore has a limit point $a^*\in \A$.
        \label{s:limit-point}
        \end{enumerate}
\end{thm}
\Pf For each $x$, the set of optimal actions $A_\alpha(x)$ in state $x$ does not
change if a constant is added to the cost function $c$. Therefore, we assume
without loss of generality that the cost function $c$ is nonnegative. Fix  $x\in\X$
and $\epsilon^*>0$. Since $x$ is fixed, we sometimes omit it.  For $\alpha\in
[0,1)$ and $a\in \A$ define
\begin{align*}
    U(x) & := \sup_{\alpha\in [0,1)} u_\alpha(x), \\
        f_\alpha(a) &:= c(x,a)+\alpha\int_\X v_\alpha(y) q(dy|x,a),\\
        g_\alpha(a) & := c(x,a)+\alpha\int_\X u_\alpha(y) q(dy|x,a).
\end{align*}
Observe that $g_\alpha(a)=f_\alpha(a)-\alpha m_\alpha$ and
\begin{align*}
    A_\alpha(x) & =\Big\{a\in\A: f_\alpha(a)=\min_{b\in\A} f_\alpha(b)\Big\}=\Big\{a\in\A:
        g_\alpha(a)=\min_{b\in\A} g_\alpha(b)\Big\}.
\end{align*}

Assumption~{\bf B} implies that $U(x)<+\infty,$ and Theorem~\ref{teor3} implies
that $\lim_{\alpha\uparrow 1} (1-\alpha)m_\alpha=w^*.$  As shown in Feinberg et al.~\cite[the
first displayed formula on p. 602]{FKZMDP}, there exists $\alpha_0\in [0,1)$ such that
for $\alpha\in [\alpha_0,1)$,
\begin{align*}
        w^*+\epsilon^*+U(x)\ge (1-\alpha)m_\alpha+u_\alpha(x) & = \min_{a\in \A}
                g_\alpha(a),
\end{align*}
This implies that for $\alpha\in [\alpha_0,1)$
\begin{align*}
        A_\alpha(x)\subseteq {\cal D}_{g_\alpha}(\lambda_1;\A)\subseteq {\cal
        D}_{g_0}(\lambda_1;\A),
\end{align*}
where the definition of the level sets ${\cal D}_\cdot(\cdot;\cdot)$ is given in
\eqref{def-D}, $\lambda_1:=w^*+\epsilon^*+U(x),$ and the second inclusion holds
because the function $u_\alpha$ takes nonnegative values.  Recall that
$f_0(a)=g_0(a)=c(x,a),$ $a\in\A,$ and the function $c(x,\cdot):\A\to\overline\R$ is
inf-compact.  Therefore, ${\cal D}_{f_0}(\lambda;\A)={\cal D}_{g_0}(\lambda;\A)$
and this set is compact for each $\lambda\in\R.$ In addition, for all $\alpha\in
[0,\alpha_0),$
\begin{align*}
    v_{\alpha_0}(x)\ge v_{\alpha}(x) & = \min_{a\in\A} f_\alpha(a),
\end{align*}
where the inequality holds because one-step costs $c$ are nonnegative. The equality is
simply the optimality equation~\eqref{eq5a}. This implies that for  $\alpha\in
[0,\alpha_0)$
\begin{align*}
        A_\alpha(x)\subseteq {\cal D}_{f_\alpha}(v_{\alpha_0}(x);\A)\subseteq {\cal
                D}_{f_0}(v_{\alpha_0}(x);\A).
\end{align*}
Let $D^*(x):= {\cal D}_{g_0}(\lambda_1;\A)\cup {\cal
D}_{f_0}(v_{\alpha_0}(x);\A).$  This set is compact as it is the union of two compact
sets, and $A_\alpha(x)\subseteq D^*(x)$ for all $\alpha\in [0,1).$  Statement
\ref{s:compact} is proved, and it implies Statement \ref{s:limit-point}. \EndPf

\section{MDPs Defined by Stochastic Equations}\label{S5}
Let $\Sp$ be a metric space, ${\cal B}(\Sp)$ be its Borel $\sigma$-field, and $\mu$
be a probability measure on   $(\Sp,{\cal B}(\Sp))$.  Consider a stochastic sequence
$\{x_\n,t \geq 0\}$ whose dynamics are defined by the stochastic equation
\begin{equation}\label{eq:stochs}
x_{\n+1}={\bf f}(x_\n,a_\n,\xi_{\n+1}),\qquad \n=0,1,\ldots,
\end{equation}
where $\{\xi_t,t \geq 1\}$ are independent and identically distributed random
variables with values in $\Sp,$ whose distributions are defined by the probability
measure $\mu,$ and ${\bf f}:\X\times\A\times\Sp\to\X$ is a continuous mapping.  This
equation defines the transition probability
\begin{equation}\label{eq:defqe}
q(B|x,a)=\int_\Sp {\bf I}\{{\bf f}(x,a,s)\in B\}\mu(ds),\qquad B\in {\cal B}(\X),
\end{equation}
from $\X\times\A$ to $ \X,$  where ${\bf I}$ is the indicator function.
%
\begin{lem}\label{l:wcq}
The transition probability $q$ is weakly continuous in $(x,a)\in \X\times\A$.
\end{lem}
\Pf
For a closed subset $B$ of $\X$ and for two sequences $x^{(n)}\to x$ and $a^{(n)}\to a$ as $n\to +\infty$ defined on $\X$ and $\A$ respectively,
\begin{align*}\limsup_{n\to\infty} q(B|x^{(n)},a^{(n)})&=\limsup_{n\to\infty}\int_\Sp {\bf I}\{{\bf f}(x^{(n)},a^{(n)},s)\in B\}\mu(ds)\\ &\le\int_\Sp \limsup_{n\to\infty}{\bf I}\{{\bf f}(x^{(n)},a^{(n)},s)\in B\}\mu(ds)\le q(B|x,a),
\end{align*}
where the first inequality follows from Fatou's lemma and the second follows from
(\ref{eq:defqe}) and upper semi-continuity of the function ${\bf I}\{{\bf f}(x,a,s)\in B\}$ on $\X\times\A\times\Sp$
for a closed set $B.$
The weak continuity of $q$ follows from Billingsley~\cite[Theorem 2.1]{Bil}
\EndPf 

\begin{cor}
    Consider an MDP $\{\X,\A,q,c \}$  with the transition function $q$ defined in
    \eqref{eq:defqe} for the continuous function $\bf f$ introduced in \eqref{eq:stochs}  and with the nonnegative $\K$-inf compact cost function $c$.  This MDP
    satisfies Assumption {\bf W*} and therefore the conclusions of
    Theorem~\ref{prop:dcoe} hold.
\end{cor}
\Pf Assumption {\bf W*}\ref{state:Wstar1} is assumed in the corollary. Assumption
{\bf W*}\ref{state:Wstar2}  holds in view of Lemma~\ref{l:wcq}. \EndPf

For inventory control problems, MDPs are usually defined by particular forms of
\eqref{eq:stochs}.  In addition, 
the cost function $c$ has the form
\begin{align}\label{e:invccost}
    c(x,a) & =C(a)+H(x,a),
\end{align}
where $C(a)$ is the ordering cost and $H(x,a)$ is either holding/backordering cost or
expected holding/back-ordering cost at the following step. For simplicity we assume
that the functions take nonnegative values. These functions are typically inf-compact.
If $C$ is lower semi-continuous and $H$ is inf-compact, then $c$ is inf-compact
because $C$ is lower semi-continuous as a function of two variables $x\in\X$ and
$a\in\A,$ and a sum of a non-negative lower semi-continuous function and an
inf-compact function is an inf-compact function.   However, as stated in the following
theorem, for discounted problems the validity of Assumption~{\bf W*} and therefore
the validity of the optimality equations, existence of optimal policies, and convergence
of value iteration take place even under weaker assumptions on the functions $C(a)$
and $H(x,a).$
\begin{thm}\label{t:cplush}
    Consider an MDP $\{\X,\A,q,c \}$  with the 
    transition function $q$ defined in
    \eqref{eq:defqe} and  cost function $c$ defined in \eqref{e:invccost}.  If either of
    the following two assumptions holds:
    \begin{enumerate}[label=(\arabic*)]
      \item the function $C:\A\to [0,\infty]$ is lower semi-continuous and the function
          $H:\X\times\A\to [0,\infty]$ is $\K$-inf-compact; \label{state:Kinf1}
      \item the function $C:\A\to [0,\infty]$ is inf-compact and the function
          $H:\X\times\A\to [0,\infty]$ is lower semi-continuous; \label{state:lsc1}
    \end{enumerate}
    \noindent  then Assumption {\bf W*} holds and therefore the conclusions of
    Theorems~\ref{prop:dcoe}\ref{state:finite1}-\ref{state:infinite}, \ref{t:lima}
    and Corollary~\ref{prop:dcoeC}\ref{cor:finite1}-\ref{cor:terminal}  hold.  Furthermore, if either of
    the following two assumptions holds:
      \begin{enumerate}[label = (\roman*)]
      \item the function $C:\A\to [0,\infty]$ is lower semi-continuous and the function  $H:\X\times\A\to [0,\infty]$ is inf-compact; \label{state:inf1}
      \item the function $C:\A\to [0,\infty]$ is inf-compact and the function
          $H^*:\A\times\X\to [0,\infty]$ is $\K$-inf-compact, where
          $H^*(a,x):=H(x,a)$ for all $x\in\X$ and all $a\in \A;$ \label{state:inf2}
      \end{enumerate}
    then the function $c$ is inf-compact and therefore  the conclusions of
    Theorems~\ref{prop:dcoe}, \ref{t:lima} and Corollary~\ref{prop:dcoeC} 
    hold.
\end{thm}
\Pf  Lemma~\ref{l:wcq} implies the weak continuity of the transition function $q.$
The definition of a $\K$-inf-compact function implies directly that the function
$C^*(x,a):=C(a)$ is  $\K$-inf-compact on $\X\times\A,$ if the function $C:\A\to
[0,\infty]$ is inf-compact.   Thus under assumptions  \ref{state:Kinf1} or
\ref{state:lsc1}, $c$ is a $\K$-inf-compact function because it is  a sum of a
nonnegative lower semi-continuous function  and a $\K$-inf-compact function. In
addition, under assumption~\ref{state:inf1}, as explained in the paragraph
preceding the formulation of the theorem, the one-step cost function $c$ is
inf-compact. Under assumption~\ref{state:inf2}, the function $c$ is inf-compact
because of the following arguments. Let  $c^*(a,x):=C(a)+H^*(a,x)$ for all
$(a,x)\in\A\times\X.$  The function $c^*:\A\times\X\to [0,\infty]$ is lower
semi-continuous if either assumption~\ref{state:Kinf1} or
assumption~\ref{state:lsc1} holds.  Since $c(x,a)=c^*(a,x)$ for all $x\in\X$ and
$a\in\A,$ the function $c:\X\times\A\to [0,\infty]$ is inf-compact if and only is the
function  $c^*:\A\times\X\to [0,\infty]$ is inf-compact.  The function $c^*$ is a
sum of the nonnegative lower semi-continuous function $C$ and the
$\K$-inf-compact function $H^*.$ Therefore, $c^*$ is $\K$-inf-compact.    Consider
an arbitrary $\lambda\in\R.$ Since $c^*(a,x)\ge C(a)>\lambda$ for $a\notin {\cal
D}_C(\lambda;\A),$ then ${\cal D}_{c^*}(\lambda;\A\times\X)={\cal
D}_{c^*}(\lambda;{\cal D}_C(\lambda;\A)\times\X),$ and this set is compact
because the set ${\cal D}_C(\lambda;\A)$ is compact and the function $c^*$ is
$\K$-inf-compact.  Thus the functions $c^*$ and $c$ are inf-compact. \EndPf

\begin{rem}
In view of Theorem~\ref{prop:dcoe}, Assumption {\bf W*} implies the existence of
optimal policies for the expected total discounted cost criterion. It is also possible
to derive sufficient conditions for the validity of Assumptions~{\bf G} and {\bf B}
and therefore for the existence of stationary optimal policies for the average costs
per unit time criterion.  However, this is more subtle than for Assumption~{\bf
W*}, and in this paper we verify Assumptions~{\bf G} and {\bf B} directly for the
periodic review inventory control problems.
\end{rem}

\section{Optimality of $(s,S)$ Policies for Setup-Cost Inventory
Control Problems}
\label{S6} In this section we consider a
discrete-time periodic-review inventory control problem with backorders, prove
the existence of an optimal $(s,S)$ policy, and establish several relevant results. For this problem, the state space is $\X:=\R,$ the action space is $\A:=\R^+,$ and the dynamics are
defined by the following stochastic equation
\begin{align}\label{e:invdinv}
x_{\n+1} &=x_\n+a_\n-D_{\n+1}, \quad \n=0,1,2,\ldots,
\end{align}
where $x_\n$ is the inventory at the end of period $\n$, $a_\n$ is the decision on how
much should be ordered, and $D_\n$ is the demand during period $\n$.  The demand
is assumed to be i.i.d.  In other words, if we change the notation $\xi_\n$ to
$D_{\n+1},$ the dynamics are defined by equation~\eqref{eq:stochs} with
${\bf f}(x,a,D)=x+a-D.$  Of course, this function is continuous.

The model has the following decision-making scenario: a
decision-maker views the current inventory of a single commodity and
makes an ordering decision. Assuming zero lead times, the products
are immediately available to meet demand. Demand is then realized,
the decision-maker views the remaining inventory, and the process
continues. Assume the unmet demand is backlogged and the cost of
inventory held or backlogged (negative inventory) is modeled as a
convex function. The demand and the order quantity are assumed to be
non-negative. 
The
dynamics of the system are defined by \eqref{e:invdinv}. Let
\begin{itemize}
  \item $\alpha \in (0,1)$ be the discount factor,
  \item $K \geq 0$ be a fixed ordering cost,
  \item $\c> 0$ be the per unit ordering cost,
  \item $D$ be a nonnegative
  random variable with the same distribution as $D_\n,$
  \item $h(\cdot)$ denote the holding/backordering cost per
  period.  It is assumed that $h:\R\to \R^+$ is a convex function, 
   $h(x) \to
  \infty$ as $|x| \to \infty,$ and $\E h(x-D) < \infty$ for all $x \in \R.$
\end{itemize}
Note that $\E D < \infty$ since, in view of Jensen's inequality,
 $h(x-\E D)\le\E h(x-D)<\infty.$ Without loss of generality, assume that $h$ is nonnegative, $h(0) = 0,$ and $h(x)>0$ for $x<0.$  Otherwise, let $x^*\in\R$ be a point, at which the function $h$ reaches its minimum value on $\R.$ Define the variable $\bar{x}:= x-x^*$ and the function ${\bar h}(\bar{x}):= h(\bar{x}+x^*)-h(x^*), $ $\bar{x}\in\R.$  Then ${\bar h}$ is a nonnegative convex function with   ${\bar h}({\bar x}) \to
  \infty$ as $|{\bar x}| \to \infty,$   ${\bar h}(0)=0,$ and ${\bar h}({\bar x})>0$ for ${\bar x}<0.$

 The cost function $c$ for this model is defined in \eqref{e:invccost} with
with $C(a):=K 1_{\{a>0\}} + \c a$ and $H(x,a):=\E h(x+a-D).$ The function
$C:\A\to \R^+$ is inf-compact. In fact, it is
continuous at $a>0$ and lower semi-continuous at $a=0.$  The function $H^*:\A\times\X\to \R^+,$ where $H^*(a,x):=H(x,a)$ for all $(a,x)\in\A\times\X,$ is
$\K$-inf-compact because of the properties of the function $h.$ Theorem~\ref{t:cplush} (case (ii)) 
implies that the function $c$ is inf-compact.  Therefore, in view of Proposition~\ref{Lemma3.2}, the function $c$ is $\K$-inf-compact.

The problem is posed with $\X = \R$ and $\A=\R^+.$ However, if the
demand and action sets are integer or
lattice, the model can be restated with $\X = \Z,$ where $\Z$ is the set of integer numbers, and $\A=\{0,1,\ldots\};$ see Remark
\ref{Rint}.

Consider the following corollary from Theorems~\ref{prop:dcoe}, \ref{t:lima}, and \ref{t:cplush}.

\begin{cor}\label{corcor}
    For the inventory control model, Assumption \textbf{W*} holds and the one-step
    cost function $c$ is inf-compact. Therefore, the conclusions of
    Theorems~\ref{prop:dcoe}, \ref{t:lima} and Corollary~\ref{prop:dcoeC} hold.  
\end{cor}
\Pf The validity of Assumption~\textbf{W*} and inf-compactness of $c$ follow from
Theorem~\ref{t:cplush} (case (ii)). \EndPf

Consider the renewal process
\begin{align}\label{renewal}
    {\bf N}(t) & := \sup\{n| \S_n \leq t\},
\end{align}
where $t\in\R^+,$ $\S_0=0$ and $\S_n= \sum_{j =1}^{n} D_j$ for $n>0.$  Observe that, if $P(D>0)>0,$ then
$\E {\bf N}(t) < \infty$ for each $ t\in\R^+;$ Resnick~\cite[Theorem
3.3.1]{res92}.  Thus, Wald's identity yields that for all $y\in\R^+$
\begin{align}\label{exp:d}
    \E {\bf S}_{{\bf N}(x)+1} & = \E ({\bf N}(y)+1) \E D < +\infty.
\end{align}
We next state a useful lemma.
\begin{lem}\label{lem:E(x)}
    For fixed initial state $x,$ if $P(D>0)>0,$ then
    \begin{align}\label{eq:defnH}
        E_y(x) & := \E h(x - \S_{{\bf N}(y)+1}) < +\infty,
    \end{align}
    where $0 \leq y < +\infty$.
\end{lem}
\Pf Define
\begin{align*}
    h^*(x) & := \begin{cases}
        h(x) & \text{for $x \leq 0$},\\
        0 & \text{for $x > 0$}.
    \end{cases}
\end{align*}
Observe that it suffices to show that
\begin{align}\label{estar}
    E_y^*(x) & := \E h^*(x - \S_{{\bf N}(y)+1}) < +\infty.
\end{align}
Indeed, for $Z=x - \S_{{\bf N}(y)+1}$,
\begin{equation*}
E_y(x)=\E 1\{Z\le 0\}h^*(Z)+ \E 1\{Z >0 \}h(Z)\le E_y^*(x)+h(x).
\end{equation*}
To show that $E_y^*(x)<+\infty,$ we shall prove that
\begin{align}\label{ineqmain}
    \E h^*(x-\S_{{\bf N}(y)+1})\le (1+\E {\bf N}(y))\E h^*(x-y-D)<+\infty.
\end{align}
Define the function $f(z)=h^*(x-y-z).$ This function
is nondecreasing and convex. Since $f$ is convex, its derivative
exists almost everywhere. Denote the excess of ${\bf N}(y)$ by
$R(y):=\S_{{\bf N}(y)+1} -y$. According to Gut~\cite[p. 59]{gut}, for $t,y\in\R^+$
\begin{align*}
    P\{R(y)>t\} & =1-F_D(y+t)+\int_0^y (1-F_D(y+t-s))dU(s),
\end{align*}
where $F_D$ is the distribution function of $D$ and $U(s)=\E {\bf N}(s)$ is the renewal function.  Thus, 
\begin{equation}\label{eqj1j2}
\E h^*(x-\S_{{\bf N}(y)+1})= \E h^*(x-y-R(y))=\E f(R(y))=\int_0^\infty
f^\prime(t) P\{R(y)>t\}dt=J_1+J_2,
\end{equation}
where $J_1=\int_0^\infty f^\prime (t)(1-F_D(y+t))dt$,
$J_2=\int_0^\infty f^\prime
(t)\left(\int_0^y(1-F_D(y+t-s))dU(s)\right)dt$, and the third equality
in \eqref{eqj1j2} holds according to Feinberg~\cite[p. 263]{fe94}. Note that
since $F_D$ is non-decreasing,
\begin{align}\label{ineqj1}
J_1 & \le  \int_0^\infty f^\prime (t)(1-F_D(t))dt=\E f(D)=\E
h^*(x-y-D)\le \E
h(x-y-D)<+\infty,
\end{align}
where the first equality follows from \cite[p. 263]{fe94}.
Similarly, by applying Fubini's theorem,
\begin{align}\label{ineqj2}
J_2 & =\int_0^y \left(\int_0^\infty f^\prime(t)(1-F_D(y+t-s))dt\right)dU(s)\nonumber \\
 & \le \int_0^y\left(\int_0^\infty f^\prime (t) (1-F_D(t))dt\right)dU(s)= \E
f(D)\E U(y)=\E h^*(x-y-D)\E {\bf N}(y).
\end{align}
Combining \eqref{eqj1j2}-\eqref{ineqj2} yields \eqref{ineqmain}. \EndPf

The following proposition is useful for the average-cost criterion.  In addition to this proposition, observe that the case $D=0$ almost surely is trivial for this criterion.  In this case, the policy $\phi,$ ordering up to the level 0, if $x<0,$ and doing nothing otherwise, is average-cost optimal.  For this policy $w(x)=w^\phi(x)= 0,$ if $x\le 0,$ and $w(x)=w^\phi(x)=h(x),$ if $x>0.$ Observe that $\phi$ is the $(0,0)$ policy.  Since $w(x)$ depends on $x,$ then Theorem~\ref{teor3} implies that Assumption ${\bf B}$ does not hold when $D=0$ almost surely.

\begin{prop}\label{prop:fixedr}  
The inventory control model satisfies Assumption~\textbf{G}  and, therefore,
 the conclusions of
    Theorem~\ref{lemC} hold.  Furthermore, if  $P(D>0)>0,$ then Assumption~\textbf{B} is satisfied and
  the conclusions of
    Theorems~\ref{teor3}, \ref{l:averoptt} and \ref{l:averopttt} hold.
\end{prop}
\Pf Consider the policy $\phi$  that orders up to the level 0, if the
inventory level is less than $0,$ and does nothing otherwise. Then
$w^{\phi}(0) =KP(D>0) +\c \E D+  \E h(-D)<+\infty.$
That is, Assumption~\textbf{G} holds.

In view of Corollary~\ref{corcor}, Theorem~\ref{prop:dcoe} implies that for every
discount factor $\alpha\in [0,1)$ there exists a stationary discount-optimal policy
$\phi^\alpha$. Theorem~\ref{lemC} implies that $\cup_{\alpha\in
[0,1)}X_\alpha\subseteq \mathcal{K}$ for some $\mathcal{K}\subseteq \R$. Let
$[x^*_L,x^*_U]$ be a bounded interval in $\R$ such that $\mathcal{K} \subseteq
[x^*_L,x^*_U]$. Thus,
\begin{align*}
    \cup_{\alpha\in [0,1)}X_\alpha \subseteq [x^*_L,x^*_U].
\end{align*}
For a discount factor $\alpha \in [0,1),$ fix a stationary optimal policy $\phi^\alpha$
and a state $x^\alpha\in [x^*_L,x^*_U]$ such that $v_\alpha(x^\alpha)=m_\alpha.$
Observe that $\phi^\alpha(x^\alpha)=0.$ Indeed, let $\phi^\alpha(x^\alpha)=a>0$.
We have
\begin{align*}
    v_\alpha(x^\alpha) & =K+ca+h(x^\alpha+a-D)+\alpha\E v_\alpha(x^\alpha+a-D) \\
        & > K+c\Big(\frac{a}{2}\Big)+h((x^\alpha+\frac{a}{2})+\frac{a}{2}-D)
            +\alpha\E v_\alpha((x^\alpha+\frac{a}{2})+\frac{a}{2}-D)\ge
            v_\alpha(x^\alpha+\frac{a}{2}),
\end{align*}
where the second inequality follows since the optimal action in state $x^{\alpha} +
\frac{\alpha}{2}$ may not be to order $\frac{a}{2}$. The inequality
$v_\alpha(x^\alpha)
>v_\alpha(x^\alpha+\frac{a}{2})$ contradicts $v_\alpha(x^\alpha)=m_\alpha$.

Let $\sigma$ be the policy defined by the following rules depending on the initial state
$x:$ (i) if $x< x^\alpha,$ then at the initial time instance $\sigma$ orders up to a
level $x^\alpha$ and then switches to $\phi^\alpha,$  and (ii)  if $x\ge x^\alpha,$
the policy $\sigma$ does not order as long as the inventory level is greater than or
equal to $x^\alpha$ and starting from the time, when the inventory level becomes
smaller than  to $x^\alpha,$ the policy $\sigma$ behaves as described in (i)  starting
from time 0.

For $x < x^\alpha,$
\begin{equation}\label{eq6.10}v_\alpha^\sigma(x)= K+\c (x^\alpha-x)+v_\alpha(x^\alpha)\le K+\c (x^*_U-x)+m_\alpha.
\end{equation}
The inequality in \eqref{eq6.10} yields for $x < x^\alpha,$
\begin{equation}\label{bound}
     v_{\alpha}(x) - m_{\alpha} \leq v_{\alpha}^{\sigma}(x) - m_{\alpha}
    \leq 
     K + \c (x^*_U-x)<+\infty.
\end{equation}
For $x \ge x^\alpha,$
\begin{equation}\label{eq6.11}
v_\alpha(x)\le v_\alpha^\sigma(x)=\E\big[\sum_{t=1}^{{\bf N}(x-x^\alpha)+1} \alpha^{t-1}h(x_t) +\alpha^{{\bf N}(x-x^\alpha)+1}[K+\c (x^\alpha-x_{{\bf N}(x-x^\alpha)+1})+v_\alpha(x^\alpha)]\big].
\end{equation}
Let $E(x):=h(x)+E_{x-x^*_L}(x)<\infty,$ where the function $E_y(x)$ is defined in
\eqref{eq:defnH} and its finiteness is stated in Lemma~\ref{lem:E(x)}. Since the
nonnegative function $h$ is convex, then  for $x_t=x-{\bf S}_t,$ $t=1,\ldots, {\bf
N}(x-x^*_L)+1,$ \begin{equation}\label{est00sum} 0\le h(x_t)\le \max\{h(x-{\bf
S}_{{\bf N}(x-x^*_L)+1}), h(x)\}\le h(x)+h(x-{\bf S}_{{\bf
N}(x-x^*_L)+1})\end{equation} and
 \begin{equation}\label{est0sum}\E h(x_t)\le h(x)+\E h(x-{\bf S}_{{\bf N}(x-x^*_L)+1})= E(x). \end{equation}
 Observe that
 \begin{equation}\label{est1sum}
 \E\big[\sum_{t=1}^{{\bf N}(x-x^\alpha)+1} \alpha^{t-1}h(x_t)\big]\le \E\big[\sum_{t=1}^{{\bf N}(x-x^*_L)+1} h(x_t)\big]\le  E(x)(1+\E {\bf N}(x-x^*_L)),
 \end{equation}
 where the first inequality follows from  $x^*_L\le x^\alpha$ and $\alpha\in [0,1);$  the second inequality follows from  $x^*_L\le x^\alpha,$ \eqref{est00sum},\eqref{est0sum}, and Wald's identity.  In addition,
 \begin{equation}\label{est2sum}\begin{aligned}
 \E[\alpha^{{\bf N}(x-x^\alpha)+1}[K+&\c (x^\alpha-x_{{\bf N}(x-x^\alpha)+1})+v_\alpha(x^\alpha)]\le K+\c(x^\alpha-x+\E{\bf S}_{{\bf N}(x-x^\alpha)+1})+m_\alpha\\
& \le K+\c(1+\E{{\bf N}(x-x^*_L)})\E D+m_\alpha,\end{aligned}
 \end{equation}
where the first inequality follows from $\alpha\in [0,1),$  $x_t=x-{\bf S}_t,$ and $v_\alpha (x^\alpha)=m_\alpha;$ the second inequality follows from $x\ge x^\alpha\ge x^*_L $ and Wald's identity. Formulae \eqref{eq6.11}, \eqref{est1sum}, and  \eqref{est2sum} imply that for $x\ge x^\alpha$
\begin{equation}\label{esecesss}
v_\alpha(x)-m_\alpha\le  K+(E(x)+\c\E D)(1+\E{{\bf N}(x-x^*_L))}<+\infty.
\end{equation}
 Inequalities \eqref{bound} and \eqref{esecesss} imply that Assumption~{\bf B} holds.
\EndPf
Consider a nonnegative, real-valued, lower semi-continuous terminal value ${\bf F}.$
%
In view of   Corollaries~\ref{prop:dcoeC}, \ref{corcor}, Theorems~\ref{prop:dcoe}, \ref{teor3}, and Proposition~\ref{prop:fixedr},   equations
\eqref{eq43333}, \eqref{eq5a} and, for the case $P(D>0)>0,$  inequality~\eqref{acoi} can be rewritten as
%
%
\begin{align}
v_{\n+1,{\bf F},\alpha}(x)& = \min \{\min_{a \ge 0} [K + G_{\n,{\bf F},\alpha}(x+a)], 
G_{\n,{\bf F},\alpha}(x)\} - \c x, \label{KDOEF2}\\
v_\alpha(x)& = \min \{\min_{a \ge 0} [K + G_{\alpha}(x+a)],
G_{\alpha}(x)\} - \c x, \label{KDOE2}\\
w + u(x) & \geq \min \{\min_{a \ge 0} [K + H(x+a)], H(x)\} - \c x, 
\label{KAOE2}
\end{align}
where  $t=0,1,\ldots$ and $w:=w(x)=w^*=\underline{w}=\overline{w},$ $x\in\X,$ and the last three equalities hold in view of \eqref{eq5.16}, and
\begin{align}\label{G111}
G_{\n,{\bf F},\alpha}(x) & := \c x + \E h(x-D) + \alpha \E v_{\n,{\bf F},\alpha}(x-D),\\
    G_{\alpha}(x) & := \c x + \E h(x-D) + \alpha \E v_{\alpha}(x-D),\label{GNG222}\\
    H(x) & := \c x + \E h(x-D) + \E u(x-D).  \label{GNG333} 
\end{align}
 We explain the correctness of \eqref{KDOEF2}.  The explanations for \eqref{KDOE2} and \eqref{KAOE2} are similar.  For this particular problem, optimality equation \eqref{eq43333}
is equivalent to $v_{\n+1,{\bf F},\alpha}(x) = \min \{\inf_{a > 0} [K +
G_{\n,{\bf F},\alpha}(x+a)], G_{\n,{\bf F},\alpha}(x)\} - \c x, $ and the internal infimum can be
replaced with the minimum in \eqref{KDOEF2} because of the following two
arguments:
\begin{enumerate}[label = (\roman*)]
  \item the function $K+G_{\n,{\bf F},\alpha}(y)$ is lower semi-continuous on
      $[x,\infty)$ and $ G_{\n,{\bf F},\alpha}(y)\to\infty$ as $y\to\infty,$ and
  \item $K+ G_{\n,{\bf F},\alpha}(x)\ge G_{\n,{\bf F},\alpha}(x)$ since $K\ge 0.$
\end{enumerate}
We remark that, in general,  while equations \eqref{KDOEF2} and \eqref{KDOE2}
are the necessary and sufficient conditions of optimality, inequality \eqref{KAOE2}
is the sufficient condition of optimality. Also, if $P(D=0)=1,$ then inequality
\eqref{KAOE2} does not hold because $w(x)$ is not a constant function, as
explained before Proposition~\ref{prop:fixedr}.
\begin{cor}\label{corollary 6.4}
    Let $\alpha\in [0,1).$  The following statements hold:
    \begin{enumerate}[label=(\alph*)]
        \item the function $G_\alpha(\cdot)$ is lower semi-continuous,
        \item if $\bf F$ is nonnegative, real-valued, and lower semi-continuous, then
            the functions $\{G_{\n,{\bf F},\alpha}(\cdot)\}_{\n=0,1,\ldots}$ are lower
            semi-continuous, and
        \item if $P(D>0)>0$, then $H$ is lower semi-continuous.
        \end{enumerate}
\end{cor}

\Pf In view of (\ref{G111})--(\ref{GNG333}), each of these functions is a sum of
several functions, two of which are continuous and the third one is lower
semi-continuous, as follows from Corollary~\ref{corcor} and from
Proposition~\ref{prop:fixedr}.\EndPf
\begin{lem}\label{l:finggn}
Let $\alpha\in [0,1).$ Then $G_\alpha(x)<+\infty$ for all $x\in\X.$ Furthermore, if
$0\le {\bf F}(x)\le v_\alpha(x)$ for all $x\in\X,$  then     $G_{\alpha,{\bf
F},\n}(x)<+\infty$ for all $x\in\X$ and for all $\n=0,1,\ldots\ .$
\end{lem}
\Pf Since $G_{\alpha,{\bf F},\n}\le G_\alpha,$ in view of \eqref{GNG222}, the lemma
follows from  $E v_\alpha(x-D)<+\infty.$  To prove this inequality, consider the policy
$\phi$  that orders up to the level 0 if the inventory level is non-positive and orders
nothing otherwise. For $x \leq 0$
\begin{align}\label{e:innegx}
    v_{\alpha}(x) \leq v^{\phi}_{\alpha}(x) \leq K - \c x +
    \frac{\alpha(K + \c \E D + \E h (-D))}{1-\alpha}.
\end{align}
Letting $B_\alpha := \frac{\alpha(K + \c \E D + \E h (-D))}{1-\alpha},$ we
have $\E v_{\alpha}(x-D) \leq K - \c\E(x-D) + B_\alpha < +\infty$.
For $x >0$,
\begin{align*}
    v_\alpha(x)  \le v_\alpha^\phi(x) & =\E\Big[\sum_{\n=1}^{{\bf N}(x)+1}\alpha^\n
        h(x-{\bf S}_\n)+ \alpha^{{\bf N}(x)+1}v_\alpha^\phi(x-{\bf S}_{{\bf N}(x)+1})\Big]\\
        & \le h(x)\E {\bf N}(x) + \E h(x-{\bf S}_{{\bf N}(x)+1})
            +K -\c (x-\E {\bf S}_{{\bf N}(x)+1})+B_\alpha<+\infty,
\end{align*}
where the second inequality follows from the facts that $\alpha^\n<1$ for $\n \geq
1$, $0\le h(x-{\bf S}_\n)\le h(x)$  for $\n=1,\ldots, {\bf N}(x),$ and
\eqref{e:innegx}. The second inequality holds because $\E {\bf N}(x)<\infty,$
Lemma~\ref{lem:E(x)}, and \eqref{exp:d}. Let $\alpha\in (0,1).$ Since
$v_\alpha^\phi(x)=\E h(x-D)+\alpha\E v_\alpha^\phi(x-D)<+\infty,$ then $\E
v_\alpha(x-D)\le\E v_\alpha^\phi(x-D)<+\infty.$  In addition, $v_0^\phi(x-D)  \le
v_\alpha^\phi(x-D)<+\infty.$ The result follows. \EndPf

Recall the following classic definition.
\begin{defn} 
For a real number $K\ge 0,$
     a function $f:\R\to\overline\R$ is called $K$-convex, if for
    each
    $x \leq y$  and for each $\lambda \in (0,1)$,
    \begin{align*}
        f((1- \lambda) x + \lambda y) \leq (1- \lambda) f(x) +
        \lambda f(y) + \lambda K.
    \end{align*}
\end{defn}
The following lemma summarizes some properties of $K$-convex functions.

\begin{lem}\label{prop:bert}
The following statements hold for a $K$-convex function $g:\R\to\R:$
\begin{enumerate}
    \item  If the function $g$ is measurable and $D$ is a random variable, then
    $\E g(x-D)$ is also $K$-convex provided $\E |g(x-D)| < \infty$ for all $x\in \R$.
    \item Suppose $g$ is inf-compact (that is,  lower semi-continuous and  $g(x)\to \infty$ as $|x| \to \infty$). Let
    \begin{align}
        S & \in \argmin_{x \in \R}\{g(x)\}, \label{defn:S}\\
        s & = \inf \{x\le S\ |\ g(x) \leq K + g(S)\}. \label{defn:s}
    \end{align}
    Then 
    \begin{enumerate}[label=(\alph*)]
        \item $g(S) + K  < g(x)$ for all $x <
        s$, \label{optS}
        \item $g(x)$ is decreasing on $(-\infty, s]$ and, therefore, $g(s)<g(x) $ for all $x<s,$ \label{decr}
        \item $g(x) \leq g(z) + K$ for all $x$ such that $s \leq x \leq z$, \label{noorder}
    \end{enumerate}

\end{enumerate}
\end{lem}
\Pf 
See
Bertsekas~\cite[Lemma~4.2.1]{bert00} and Simchi-Levi et al.~\cite[Lemma 8.3.2]{SCB} for the case of a continuous function $g.$  The proofs there with minor adjustments cover the case when $g$ satisfies the measurability and continuity properties stated in the lemma.  \EndPf
%
%
%

Consider the discounted cost problem and suppose $G_{\alpha}$ is $K$-convex, lower
semi-continuous and approaches infinity as $|x| \to \infty$. If we define $S_{\alpha}$
and $s_{\alpha}$ by \eqref{defn:S} and \eqref{defn:s} with $g$ replaced by
$G_{\alpha}$, Statement 2\ref{optS} of Lemma~\ref{prop:bert}, along with the
optimality equation~\eqref{KDOE2}, imply that it is optimal to order up to
$S_{\alpha}$ when $x < s_{\alpha}$. Statement 2\ref{noorder} of
Lemma~\ref{prop:bert} imply that it is optimal not to order when $x\ge s_{\alpha}.$
Our next goal is the established these properties of the function $G_\alpha$ and of
some relevant functions.

For a fixed ordering cost $K\ge 0$ we sometimes write $v_\alpha^K,$ $v_{t,\alpha}^K,$ $v_{t,{\bf F},\alpha}^K,$ $G^K_\alpha,$ and   $G_{t,{\bf F},\alpha}^K,$
instead of   $v_\alpha,$  $v_{t,\alpha},$ $v_{t,{\bf F},\alpha},$ $G_\alpha,$ and $G_{t,{\bf F},\alpha},$ respectively. Consider the terminal value ${\bf F}(x)=v^0_\alpha(x),$ $x\in\X.$ According to Theorem~\ref{prop:dcoe}\ref{state:fl1} and Corollary~\ref{prop:dcoeC}\ref{cor:inf-compact}, the functions $v_\alpha,$ $v^0_\alpha,$  $v_{t,\alpha},$ and $v_{t,v^0_\alpha,\alpha},$ $t=1,2,\ldots,$ are inf-compact.  

\begin{lem}\label{Lemmmamm} The following statements hold:
\begin{enumerate}[label = (\roman*)]
    \item the functions $v_\alpha$ and $v_{t,v^0_\alpha,\alpha},$ $t=0,1,\ldots,$
        are inf-compact;\label{state:v-inf}
    \item the functions  $G_\alpha$ and $G_{t,v^0_\alpha,\alpha},$
        $t=0,1,\ldots,$ are lower semi-continuous, and \label{state:g-inf}
        \begin{align*}
            \lim_{x\to +\infty} G_\alpha(x) & =  \lim_{x\to +\infty}
                G_{t,v^0_\alpha,\alpha}(x)=+\infty, \qquad t=0,1,\ldots;
        \end{align*}
    \item there exists $\alpha^*\in [0,1)$ such that $G^0_\alpha(x)\to\infty$ as
        $x\to -\infty$ for all $\alpha\in[\alpha^*,1) ;$\label{state:gfinite}
    \item for $\alpha\in [\alpha^*,1),$ where $\alpha^*$ is the constant
        $\alpha^*\in [0,1)$ whose existence is stated in
        Statement~\ref{state:gfinite}, the functions $G_\alpha(x)$  and
        $G_{\n,v_\alpha^0,\alpha}(x), $ $\n=0,1,\ldots,$ are $K$-convex and tend
        to $+\infty$ as $x\to-\infty,$  and therefore, in view of Statement
        \ref{state:g-inf}, these functions are inf-compact. Furthermore, the functions
        $v_\alpha$ and $v_{\n,v_\alpha^0,\alpha}(x), $ $\n=0,1,\ldots,$ are
        $K$-convex.  \label{state:ave-ginf}
\end{enumerate}
\end{lem}
\Pf In view of Corollary~\ref{corcor}, Statement~\ref{state:v-inf} follows from
Theorem~\ref{prop:dcoe}\ref{state:fl1} and
Corollary~\ref{prop:dcoeC}\ref{cor:inf-compact}. Statement~\ref{state:g-inf}
follows from Statement~\ref{state:v-inf}, nonnegativity of costs, and
definitions~\eqref{G111} and \eqref{GNG222}.

To prove Statement~\ref{state:gfinite} note that it is well-known that the function
$G^0_\alpha$ is convex, where $\alpha\in [0,1).$  Indeed, the function
$v^0_{0,\alpha}=0$ is convex.  For $K=0,$ equations \eqref{KDOEF2},
\eqref{G111} and induction based on Heyman and Sobel~\cite[Proposition
B-4]{heysol} imply that the functions $v^0_{\n,\alpha},$ $\n=1,2,\ldots,$ are
convex. Convergence of value iterations, stated in
Theorem~\ref{prop:dcoe}\ref{state:finite1}, implies the convexity of the functions
$v^0_\alpha.$  The convexity of $G^0_\alpha$ follows from \eqref{GNG222}.

We show by contradiction that there exists $\alpha^*\in [0,1)$
such that $G^0_\alpha$ is decreasing on an interval
$(-\infty,M_\alpha]$ for some $M_\alpha>-\infty$ when $\alpha\in
[\alpha^*,1)$. Suppose this is not the case. For $K=0$,
\eqref{KDOE2} can be written as
\begin{equation}
v_{\alpha}^0(x) =  \inf_{a \ge 0} \{ G^0_{\alpha}(x+a)\} - \c x.
\label{KDOE34}
\end{equation}
If a constant $M_\alpha$ does not exist for some $\alpha \in (0,1),$ then the
convexity  of $G^0_\alpha(x)$ implies that the policy $\psi,$ that never orders, is
optimal for the discount factor $\alpha$.  If there is no $\alpha^*$ with the
described property, Corollary~\ref{l:averopt} implies that the policy $\psi$ is
average-cost optimal.  This is impossible because, if $x$ is small enough that the
convex function $h(x)$ is decreasing at $x,$ then  $w^\psi(x)\ge \ E h(x-D)>h(x)\to
+\infty$ as $x\to-\infty,$ but, in view of Theorem~\ref{teor3}, $w(x)$ is a finite
constant.  This contradiction implies that for $\alpha\in [\alpha^*,1)$ the functions
$G^0_\alpha$ decreases when $x\in (-\infty,M_\alpha],$ where $M_\alpha>-\infty.$
The convexity of $G^0_\alpha$ implies that $G^0_\alpha(x)\to\infty$ as $x\to
-\infty.$

Let us prove Statement~\ref{state:ave-ginf}. The convergence of the functions to $+\infty,$ as $x\to -\infty,$ follows from Statement~\ref{state:gfinite} and the inequalities $G_\alpha^K(x)\ge G^0_\alpha (x)$ and $G_{t,v^0_\alpha,\alpha}^K(x)\ge G^0_\alpha (x),$ which hold for all $x\in\X.$ Indeed, the first inequality follows from $v_\alpha^K(x)\ge v^0_\alpha (x),$ $x\in\X,$ and \eqref{GNG222}.  The second inequality follows from
$v^K_{t,v^0_\alpha,\alpha}(x) \ge v^0_{t,v^0_\alpha,\alpha}(x)= v^0_\alpha (x),$ $x\in\X$, and \eqref{G111}.

Now let $\alpha\in [\alpha^*,1).$
As explained in the proof of \ref{state:gfinite},
the function  $G^0_\alpha$ is
convex and therefore it is $K$-convex. Formulae \eqref{KDOEF2}, \eqref{G111},
Heyman and Sobel \cite[Lemma 7-2, p. 312]{heysol}, and induction arguments imply
that the functions $G_{\n,v_\alpha^0,\alpha}$ and $v_{\n+1,v_\alpha^0,\alpha},$
$t=1,2,\ldots$ are $K$-convex.  In addition, $v_{\n,v_\alpha^0,\alpha}(x)\uparrow
v_{\alpha}(x)$ 
as $t\to\infty$ in view of Corollary~\ref{prop:dcoeC}\ref{cor:terminal} and since all
the costs are nonnegative. Formulae \eqref{G111}, \eqref{GNG222} and the
monotone convergence theorem imply that $G_{\n,v_\alpha^0,\alpha}(x)\uparrow
G_{\alpha}(x)$ as $t\to\infty.$ Thus, the functions $v_\alpha$ and $G_\alpha$ are
$K$-convex.
\EndPf

\begin{defn}
Let $s_\n$ and $S_\n$ be real numbers such that $s_\n\le S_\n$,
$\n=0,1,\ldots\ .$ Suppose $x_\n$ denotes the current inventory
level at decision epoch $\n$. A policy is called an $(s_\n,S_\n)$
policy at step $\n$ if it orders up to the level $S_\n$ if $x_\n <
s_\n$ and does not order when $x_\n \ge s_\n.$  A Markov 
policy is called an $(s_\n,S_\n)$ policy if it is an $(s_\n,S_\n)$
policy at all steps $\n=0,1,\ldots\ .$ 
A policy is called
an $(s,S)$ policy if it is stationary and it is an $(s,S)$ policy
at all steps $\n=0,1,\ldots\ .$
\end{defn}
The following theorem is the main result of this section.
\begin{thm} \label{th:sSdisc}
Consider $\alpha^*\in [0,1)$  whose existence is stated in
Lemma~\ref{Lemmmamm}. The following statements hold for the inventory
control problem.
    \begin{enumerate}[label= (\roman*)]
      \item For $\alpha \in [\alpha^*,1)$ and $\n=0,1,\ldots,$ define
          $g(x):=G_{\n,v_\alpha^0,\alpha}(x),$ $x\in\R$. Consider real numbers
          $S^*_{\n,\alpha}$ satisfying  \eqref{defn:S} and $s^*_{\n,\alpha}$
          defined in   \eqref{defn:s}. For each $N=1,2,\ldots,$  the $(s_t,S_t)$
          policy with $s_t=s^*_{N-\n-1,\alpha}$ and $S_t=S^*_{N-\n-1,\alpha},$
          $\n=0,1,\ldots,N-1,$ is optimal for the $N$-horizon problem with the
          terminal values ${\bf F}(x)=v^0_\alpha(x),$ $x\in\R$. \label{state:inv1}
      \item For the infinite-horizon expected total discounted cost criterion with a
          discount factor $\alpha\in [\alpha^*,1),$ define $g(x):=G_\alpha(x),$
          $x\in\R$. Consider real numbers $S_{\alpha}$ satisfying \eqref{defn:S}
          and $s_{\alpha}$ defined in   \eqref{defn:s}. The $(s_\alpha,S_\alpha)$
          policy is optimal for the discount factor $\alpha.$ Furthermore, a sequence
          of pairs $\{(s^*_{\n,\alpha},S^*_{\n,\alpha})\}_{\n=0,1,\ldots} $ is
          bounded, where $s^*_{\n,\alpha}$ and $S^*_{\n,\alpha}$ are described
          in Statement \ref{state:inv1}, $\n=0,1,\ldots\ .$ If
          $(s^*_\alpha,S^*_\alpha)$ is a limit  point of this sequence, then the
          $(s^*_\alpha,S^*_\alpha)$ policy is optimal for the infinite-horizon
          problem with  the discount factor $\alpha.$ \label{state:inv2}
      \item Consider the infinite-horizon average cost  criterion. For each $\alpha\in
          [\alpha^*,1)$, consider an optimal $(s^\prime_\alpha, S^\prime_\alpha)$
          policy for the discounted cost criterion with the discount factor $\alpha,$
          whose existence follows from Statement (ii).   Let $\alpha_n\uparrow 1,$
          $n=1,2,\ldots,$ with $\alpha_1\ge\alpha^*.$  Every sequence
          $\{(s^\prime_{\alpha_n}, S^\prime_{\alpha_n}), n \geq 1\}$ is bounded
          and  each its limit point $(s,S)$ defines an average-cost optimal $(s,S)$
          policy. Furthermore, if $P(D>0)>0,$ this policy satisfies the optimality
          inequality \eqref{KAOE2} with $u=\tilde u,$ where the function $\tilde u$ is
          defined in \eqref{eq1311} for an arbitrary subsequence
         $\{\alpha_{n_k}\}_{k=1,2,\ldots}$ of $\{\alpha_{n},n\geq 1\}$ satisfying
         $(s,S)=\lim_{k\to\infty} (s^\prime_{\alpha_{n_k}},
         S^\prime_{\alpha_{n_k}}).$
          \label{state:inv3}
    \end{enumerate}
\end{thm}
\Pf To prove Statements \ref{state:inv1} and \ref{state:inv2}, let $\alpha\in
[\alpha^*,1).$  In view of Lemma~\ref{Lemmmamm}\ref{state:ave-ginf}, the
functions $G_\alpha$ and $G_{\n,v^0_\alpha,\alpha},$ $t=0,1,\ldots,$ are
$K$-convex and inf-compact.  The optimality of $(s_t,S_t)$ policies and $(s,S)$
policies with $s=s_\alpha$ and $S=S_\alpha$ stated in \ref{state:inv1} and
\ref{state:inv2} follows from optimality equations \eqref{KDOEF2}, \eqref{KDOE2},
Lemma~\ref{prop:bert} with $g=G_{N,v^0_\alpha,\alpha}$ and $g=G_\alpha$
respectively, and Theorem~\ref{prop:dcoe}.

Consider now the remaining claims in \ref{state:inv2}.  Since $G^0_\alpha(x)\le
G_{\n,v_\alpha^0,\alpha}(x)\le G_{\n+1,v_\alpha^0,\alpha}(x)\le G_{\alpha}(x),$
$x\in\R,$ the points $s^*_{\n,\alpha}$ and $S^*_{\n,\alpha}$ belong to the compact
set $\{x\in\R: G_\alpha^0(x)\le K+\min_{x\in\R}G_\alpha(x) \}.$ Therefore, the
sequence $\{(s^*_{\n,\alpha},S^*_{\n,\alpha})\}_{t=0,1,\ldots}$ is bounded and
has a limit point $(s^*_\alpha,S^*_\alpha).$ The function ${\bf
F}(x)=v_\alpha^0(x)$ satisfies inequalities in \eqref{eq:spdouco}, and therefore the
assumptions of Theorem~\ref{t:lima} hold. Theorem~\ref{t:lima} implies that, for the
infinite-horizon problem with the discount factor $\alpha,$  the following decisions are
optimal for the corresponding states: no inventory should be ordered for
$x>s^*_\alpha$ and the inventory up to the level $S^*_\alpha$ should be ordered for
$x<s_\alpha^*.$  This implies that $G_\alpha(x)\le K+ G_\alpha(S^*_\alpha)$ for
$x\in (s^*_\alpha, S^*\alpha).$  Lower semi-continuity of $G_\alpha(x)$ implies that
$G_\alpha(s^*_\alpha)\le K+ G_\alpha(S^*_\alpha).$ Thus, the decision, that
inventory should not be ordered, is optimal at $x=s^*_\alpha.$  That is, the
$(s^*_\alpha,S^*_\alpha)$ policy is optimal for the infinite-horizon problem with the
discount factor $\alpha.$

It remains to prove Statement \ref{state:inv3}. Let $P(D>0)>0.$ We start with the
proof that the sequence $\{(s^\prime_{\alpha_n},
S^\prime_{\alpha_n})\}_{n=1,2,\ldots}$ is bounded. First, we prove that the
sequence $ \{s^\prime_{\alpha_n}, n \geq 1\}$ is bounded below. If this is not true,
then $\lim_{k\to\infty} s^\prime_{\alpha_{n_k}}=-\infty$ for some $n_k\to\infty$ as
$k\to \infty.$ This means that for each $x\in\R$ there is a natural number $k(x)$ such
that $x>s^\prime_{\alpha_{n_k}}$ for $k\ge k(x).$ Therefore,
$0\in\A_{\alpha_{n_k}}(y),$ $k\ge k(x), $ for all $y\ge x.$
Corollary~\ref{l:averopt}\ref{state:averopt2} implies that the action $0\in
A^*_{\tilde u}(y)$ for all $y>x,$ where $\tilde{u}$ is defined in \eqref{eq1311} for
the sequence of discount factors $\{\alpha_{n_k}, k\geq 1\}$.  Since $x\in\R$ is
arbitrary,   $0\in A^*_{\tilde u}(y)$ for all $y\in\R.$   This means that the policy
$\psi,$ that never  orders inventory, is optimal for average costs per unit time.
However,
\begin{align*}
    w^\psi(x)\ge \E h(x-{\bf S}_n) & \ge h(x-n\E D).
\end{align*}
Letting $n\to\infty$ on the right hand side yields $w^\psi(x)=+\infty$ for all $x\in\R$.
In view of Assumption {\bf G}, that holds for the inventory control problem,
$w^\psi(x)<+\infty$ for some $x\in\R.$ Thus  the sequence $\{s^\prime_{\alpha_n},
n \geq 1 \}$ is bounded.

Second, we prove that the sequence $ \{S^\prime_{\alpha_n}, n \geq 1\}$ is also
bounded.  Let $x\in\R$ be a lower bound for $\{s^\prime_{\alpha_n},  n \geq 1\}.$
Thus, $a^{(n)}:= (S^\prime_{\alpha_n} -x)\in A_{\alpha_n}(x).$ In view of
Theorem~\ref{l:averopttt},  the sequence $\{a^{(n)}, n \geq 1\}$ is bounded. This
implies that the sequence $\{S^\prime_{\alpha_n}, n \geq 1\}$ is bounded as well.

Consider a subsequence $\alpha_{n_k}\uparrow1$ such that
$(s^\prime_{\alpha_{n_k}},S^\prime_{\alpha_{n_k}})\to (s^\prime, S^\prime)$
as $k\to\infty.$ Corollary~\ref{l:averopt}\ref{state:averopt2} implies that $0\in
A^*_{\tilde u}(x),$ if $x>s^\prime,$ and $S^\prime-x\in A^*_{\tilde u}(x),$ if
$x<s^\prime,$ where the function $\tilde u$ is defined in \eqref{eq1311} for the
sequence  of discount factors $\{\alpha_{n_k},k \geq 1\}$. The last step is to
prove that $0\in A^*_{\tilde u}(s^\prime).$ To do this, consider a subsequence
$\{\alpha^*_n, n \geq 1\}$ such that $\alpha^*_n \to 1$ of the sequence
$\{\alpha_{n_k}, k \geq 1\}$ and a sequence $\{x^{(n)}, n \geq 1\}$ with
$x^{(n)} \to s^\prime$ such that $\tilde{u}(s^\prime)=\lim_{n\to\infty}
u_{\alpha^*_n}(x^{(n)})$.

First, consider the case when there is a sequence $\ell_k\to\infty$ such that
$x^{(\ell_k)}\ge s^\prime_{\alpha^*_{\ell_k}}$ for all $k=1,2,\ldots\ .$ In this
case, $0\in A_{\alpha^*_{\ell_k}}(x^{(\ell_k)}),$ and
Corollary~\ref{l:averopt}\ref{state:averopt2} implies that $0\in
A^*_{\tilde{u}^*}(s^\prime),$ where the function ${\tilde u}^*$ is defined in
\eqref{eq1311} for the sequence  of discount factors
$\{\alpha^*_{n_k}\}_{k=1,2,\ldots}.$  Observe that ${\tilde u}^*(s^\prime)=
{\tilde u}(s^\prime)$ and ${\tilde u}^*(x)\ge {\tilde u}(x)$ for all $x\in\R.$ This
implies $A^*_{\tilde{u}^*}(s^\prime)\subseteq A^*_{\tilde{u}}(s^\prime).$
Thus $0\in A^*_{\tilde{u}}(s^\prime).$

Second, consider the complimentary case, when there exists a number $N$ such
that $x^{(n)}<s^\prime_{\alpha^*_{n}} $ for $n\ge N.$ Let $n\ge N.$ In view of
Statement~2\ref{decr} of  Lemma~\ref{prop:bert},
$G_{\alpha^*_{n}}(x^{(n)})\ge G_{\alpha^*_{n}}(s^\prime_{\alpha^*_{n}})$.
Therefore,
 \begin{align*}
    u_{\alpha^*_{n}}(x^{(n)}) & = v_{\alpha^*_{n}}(x^{(n)})-m_{\alpha^*_{n}}
        =K+G_{\alpha^*_{n}}(S^\prime_{\alpha^*_{n}})-\c x^{(n)}
            - m_{\alpha^*_{n}}
         \ge G_{\alpha^*_{n}}(s^\prime_{\alpha^*_{n}}) -\c x^{(n)}
            - m_{\alpha^*_{n}} \\
        & \ge  v_{\alpha^*_{n}}(s^\prime_{\alpha^*_{n}})
            +\c s^\prime_{\alpha^*_{n}}-\c x^{(n)} - m_{\alpha^*_{n}}
         = u_{\alpha^*_{n}}(s^\prime_{\alpha^*_{n}})
            +\c (s^\prime_{\alpha^*_{n}}- x^{(n)}),
 \end{align*}
where the first and the last equalities follow from the definition of the functions
$u_\alpha,$ the second equality follows from \eqref{KDOE2} and from the
optimality of the $(s^\prime_{\alpha^*_{n}},S^\prime_{\alpha^*_{n}})$
policies for discount factors $\alpha^*_{n},$  the first inequality follows from Statement~2\ref{noorder} of
Lemma~\ref{prop:bert}, and the last inequality follows from \eqref{KDOE2}.
Since $s^\prime_{\alpha^*_{n}}\to s^\prime$ and  $x^{(n)}\to s^\prime,$
\begin{align*}
    \tilde{u}(s^\prime)& =\lim_{n \to\infty }u_{\alpha^*_{n}}(x^{(n)})
        =\lim_{n\to\infty } u_{\alpha^*_{n}}( s^\prime_{\alpha^*_{n}}).
\end{align*}
Moreover, since $0\in A_{\alpha^*_{n}}(s^\prime_{\alpha^*_{n}})$ for all
$n=1,2,\ldots,$  Theorem~\ref{l:averoptt}\ref{state:converge2} implies that $0\in
A^*_{\tilde u}(s^\prime).$ Thus, the $(s^\prime,S^\prime)$ policy is average-cost
optimal.

Now let $D=0$ almost surely. As explained in the paragraph preceding
Proposition~\ref{prop:fixedr}, the $(0,0)$ policy $\phi$ is average-cost optimal. Let
us prove that
\begin{equation}\label{eqs0S0}
\lim_{\alpha\uparrow 1} s_\alpha = \lim_{\alpha\uparrow 1} S_\alpha =0.
\end{equation}
Let $\alpha\in (0,1).$ Consider an arbitrary policy $\sigma.$ Since $v^\sigma(x)\ge
\frac{h(x)}{1-\alpha}= v^\phi_\alpha(x),$ when $x\ge 0,$ then
$v_\alpha(x)=\frac{h(x)}{1-\alpha}$ for all $x\ge 0.$ This formula and
\eqref{GNG222} imply $G_\alpha(x) =\c x+h(x)/(1-\alpha)$ for $x\ge 0.$  Thus, the
function $G_\alpha(x)$ is increasing, when $x\in [0,\infty).$  This implies $S_\alpha\le
0.$ Since $s_\alpha\le S_\alpha,$ then $s^*=\liminf_{\alpha\uparrow 1} s_\alpha\le
0.$  To complete the proof of \eqref{eqs0S0}, we need to show that $s^*=0.$ Indeed,
let us assume that $s^*<0.$  Fix an arbitrary $x\in (s^*,0).$ Then there exists a
sequence $\alpha_n\uparrow 0$ such that $s_{\alpha_n}\to s^*$ as $n\to\infty$ and
$s_{\alpha_n}< x,$ $n=1,2,\ldots\ .$  The $(s_{\alpha_n},S_{\alpha_n})$ policy
$\phi^n$ is optimal for the discount factor $\alpha_n,$ and this policy does not order
at the state $x,$ $n=1,2,\ldots\ .$ Therefore $v_{\alpha_n}^{\phi^n}(x)=
h(x)/(1-\alpha_n)\to +\infty$ as $n\to\infty.$ However,  $v^\phi_{\alpha_n}(x)=K-\c
x.$ This implies that the $(s_{\alpha_n},S_{\alpha_n})$ policy $\phi^n$ cannot be
optimal for a discount factor $\alpha_n> (K-\c x)/(K-\c x -h(x)).$ \EndPf

For $N=1,2,\ldots,$ we shall write $G_{N,\alpha}$ instead of
$G_{N,{\bf F},\alpha}$ if ${\bf F}(x)=0$ for all $x\in\R.$

\begin{lem}\label{lem:galpha}
       Suppose there exist $z,y\in\R$ such that $z<y$ and
            \begin{align}
                \frac{h(y) - h(z)}{y-z} < -\c. \label{eq:slopeS}
            \end{align}
        Then $G_\alpha(x)\to +\infty \text{ and } G_{N,\alpha}(x) \to +\infty$ as $|x| \to
        \infty$ for all $\alpha\in [0,1)$ and for all $N \geq 0$, and these functions are $K$-convex.
\end{lem}

\Pf 
 Observe that the assumption in Lemma~\ref{lem:galpha} is equivalent to the existence of $z,y\in\R$ such that $z<y$ and
  \begin{align}
                \frac{\E [h(y-D) - h(z-D)]}{y-z} < -\c. \label{eq:slope}
            \end{align}
Indeed, since  $h$ is convex, $h(y-d)-h(z-d)\le h(y)-h(z),$ and \eqref{eq:slopeS} implies \eqref{eq:slope}.  Also, \eqref{eq:slope}
implies that for some $d\ge 0$ inequality \eqref{eq:slopeS} holds for $y:=y-d$ and $z:=z-d.$

According to \eqref{G111}, $G_{N,\alpha}(x) \to \infty$ as $x \to\infty$ for all $N=0,1,\ldots\ .$
We show that the result continues to hold when $x \to -\infty$.
Suppose $z < y$ satisfy \eqref{eq:slope}.  Inequality \eqref{eq:slope} can be rewritten as
\[\c y  + \E h(y-D) < \c z + \E h(z-D).\]
Thus, $G_{0,\alpha}(z) > G_{0,\alpha}(y)$. Since $G_{0,\alpha}$ is
convex, then $G_{0,\alpha}(x) \to \infty$ as $x \to -\infty.$ According to \eqref{G111},
\begin{align*}
G_{N,\alpha}(x)  &= G_{0,\alpha}(x) + \alpha\E v_{N,\alpha}(x - D)\ge G_{0,\alpha}(x),\qquad N=1,2,\ldots,\\
G_{\alpha}(x)  &= G_{\alpha}(x) + \alpha\E v_{\alpha}(x - D)\ge G_{0,\alpha}(x),
\end{align*}
where $ G_{0,\alpha}(x)\to +\infty$ as $x\to -\infty.$ \EndPf

\begin{thm}\label{th:sSdisc1} Under the condition stated in Lemma~\ref{lem:galpha},  the following statements hold for each discount factor $\alpha \in [0,1)$:
\begin{enumerate}[label = (\roman*)]
  \item For $\n=0,1,\ldots$ consider real numbers
      $S_{\n,\alpha}$ satisfying      \eqref{defn:S} and   $s_{\n,\alpha}$ defined in
      \eqref{defn:s}  with $g(x)=G_{\n,\alpha}(x),$ $x\in\R.$ Then for
      every $N=1,2,\ldots$  the  $(s_t,S_t)$ policy with $s_t=s_{N-\n - 1,\alpha}$ and $S_t=S_{N-\n - 1,\alpha},$
      $\n=0,1,\ldots,N-1,$ is  optimal for the $N$-horizon problem with the
      zero terminal values.
  \item  Consider real numbers $S_{\alpha}$ satisfying \eqref{defn:S}
          and $s_{\alpha}$ defined in   \eqref{defn:s} for $g(x):=G_\alpha(x),$
          $x\in\R.$ Then the
          $(s_\alpha,S_\alpha)$ policy is optimal for the infinite-horizon problem with the discount factor $\alpha.$ 
 Furthermore, a sequence of pairs
          $\{(s_{\n,\alpha},S_{\n,\alpha})\}_{\n=0,1,\ldots} $ considered in statement (i) is bounded, and, if
          $(s^*_\alpha,S^*_\alpha)$  is a limit  point of this sequence, then the
          $(s^*_\alpha,S^*_\alpha)$ policy is optimal for the infinite-horizon problem with the discount factor
          $\alpha.$
\end{enumerate}
\end{thm}
\Pf
Observe that $G_{0,\alpha}(x)=\c x +\E h(x-D).$ This function is convex and, in view of Lemma~\ref{lem:galpha}, $G_{0,\alpha}(x)\to\infty$ as $|x|\to\infty.$
The rest of the proof coincides with the proof of Theorem~\ref{th:sSdisc} with the functions $G_{\n,v_\alpha^0,\alpha}$ replaced with the functions   $G_{\n,\alpha}.$
\EndPf

By using the results of this section, Feinberg and Liang~\cite{Fli1, Fli} obtained
additional results for the inventory control problem.   Feinberg and Liang~\cite{Fli}
described the structure of optimal policies for all values of the discount factor
$\alpha\ge 0$ for finite-horizon problems and  for all values of $\alpha\in [0,1) $ for
infinite-horizon problems.  In particular, the smallest possible values of the
discount factor $\alpha^*$ mentioned in Theorem~\ref{th:sSdisc} are computed in
\cite{Fli}.  Though the general theory of MDPs implies that the value functions
$v_{t,\alpha}(x),$ $G_{t,\alpha}(x),$  $v_\alpha(x),$ and $G_\alpha(x)$ are lower
semi-continuous in $x,$  it is proved in  \cite{Fli} that these functions are
continuous.  In particular, these continuity properties imply that, for total
discounted cost criteria with finite and infinite horizons, the decisions to order up to
the level $S$ ($S_t$) are also optimal at the states $s$ ($s_t$). Feinberg and
Liang~\cite{Fli1} proved that for the inventory control problem the average-cost
optimality inequality in \eqref{KAOE2} holds in the stronger form of the optimality
equation, the convergences $u_\alpha(\cdot)\to u(\cdot)$ and $G_\alpha(\cdot)\to
H(\cdot)$ take place, as $\alpha\uparrow 1,$ and the functions $u(x)$ and $G(x)$
are $K$-convex and continuous.  Therefore, average-cost optimal $(s,S)$ policies
can be derived from the optimality equation, and the decision to place an order up
to the level $S$ at the state $s$ is also optimal for the average-cost criterion.


\begin{rem}
This remark comments on the assumptions $\alpha\in [0,1),$ $K\ge 0,$ and $c>0.$
All the results of this paper stated for the finite horizon  hold with the same proofs
for arbitrary $\alpha \ge 0;$ see Feinberg and Liang~\cite{Fli} for detail. If $K=0,$
then it is well-known that it is possible to set $s=S$ and $s_t=S_t$ for the
corresponding optimal $(s,S)$ policies, see e.g., Heyman and
Sobel~\cite[Proposition 3-1]{heysol}, and such policies are called base stock or
$S$-policies.  Indeed, this follows from Lemma~\ref{Lemmmamm} and
\eqref{KDOEF2}, \eqref{KDOE2} for discounted problems, and then from
Theorem~\ref{th:sSdisc}\ref{state:inv3} for problems with average costs per unit
time. If $c=0,$ then Assumption~{\bf W*} holds.  In particular, the function
$c(x,a)=K1_{a>0}+\E h(x+a-D)$ is $\K$-inf-compact as a sum of a
lower-semicontinuous function and a $\K$-inf-compact function; see
Theorem~\ref{t:cplush}\ref{state:Kinf1}.  All the results formulated in the paper
for a fixed discount factor $\alpha\in [0,1)$ remain correct for $\bar c=0.$
Furthermore, inequality~\eqref{eq:slopeS} holds and therefore the conclusions of
Theorem~\ref{th:sSdisc1} hold.  However, the function $c$ is not inf-compact when
$\bar c =0.$ For example, $c(-a,a)=K+\E h(-D)\nrightarrow +\infty$ as $a\to
+\infty.$  The proof of  Assumption~{\bf B} in Proposition~\ref{prop:fixedr} is
based on Theorem~\ref{lemC}, which uses the assumption that the function $c$ is
inf-compact.  So, for the long-term average-cost criterion, the results of this paper
do not cover the case ${\bar c}=0.$
\end{rem}

\begin{rem}\label{Rint} For the inventory control problem, we
have considered an MDP with $\X = \R$ and $A(x) = \R^+$ for each $x \in \X$.
However, if the demand takes only integer values, for many problems it is natural to
consider $\X = \Z$ and $A(x) = \Z^+,$ where $\Z$ is the set of integers and
$\Z^+$ is the set of nonnegative integers.  Therefore, if the demand is integer, we
have two MDPs for the inventory control problems: an MDP with $\X = \R$ and an
MDP with $\X = \Z$. All of the results of this paper also hold for the second
representation, when the state space is integer, with a minor modification that the
action sets are integer as well. In fact the case $\X=\Z$ is slightly easier because
every function is continuous on it and therefore it is lower semi-continuous.
\end{rem}

\noindent {\bf Acknowledgement.}  Research of the first author was partially
supported by NSF grants CMMI-1335296 and CMMI-1636193.
Theorem 6.10, 6.12 of this paper and some relevant preceding constructions and
definitions are reprinted with some minor differences by permission from \cite{fe16}, where
these theorems are presented without proofs. Copyright 2016, the Institute for
Operations Research and the Management Sciences, 5521 Research Park Drive,
Suite 200, Catonsville, MD 21228 USA.

\end{document}